\documentclass[twocolumn]{mmaauth}
\usepackage{moreverb}
\usepackage{bm}

\usepackage[dvips,colorlinks,bookmarksopen,bookmarksnumbered,citecolor=red,urlcolor=red]{hyperref}

\newcommand\BibTeX{{\rmfamily B\kern-.05em \textsc{i\kern-.025em b}\kern-.08em
T\kern-.1667em\lower.7ex\hbox{E}\kern-.125emX}}

\def\volumeyear{2017}
\newtheorem{theorem}{Theorem}[section]
\newtheorem{lemma}[theorem]{Lemma}
\newtheorem{proposition}[theorem]{Proposition}
\newtheorem{corollary}[theorem]{Corollary}
\newtheorem{remark}[theorem]{Remark}

\theoremstyle{definition}
\newtheorem{definition}[theorem]{Definition}

\def\supp{\text{supp}}

\numberwithin{equation}{section}
\begin{document}

\runninghead{Ningning Han, Shidong Li, Zhanjie Song, Hong Wang}

\title{The finite steps of convergence of the fast thresholding algorithms with feedbacks}

\author{Ningning Han\affil{a},
Shidong Li\affil{b}\corrauth\, Zhanjie Song\affil{a}, Hong Wang\affil{a}}

\address{\affil{a}School of mathematical sciences, Tianjin University, Tianjin 300350, China.\affil{b}Department of Mathematics, San Francisco State University, San Francisco, CA94132,USA}

\corraddr{Department of Mathematics, San Francisco State University, San Francisco, CA94132,USA. E-mail: shidong@sfsu.edu.}

\begin{abstract}
Iterative algorithms based on thresholding, feedback and null space tuning (NST+HT+FB) for sparse signal recovery are exceedingly effective and fast, particularly for large scale problems. The core algorithm is shown to converge in finitely many steps under a (preconditioned) restricted isometry condition. In this paper, we present a new perspective to analyze the algorithm, which turns out that the efficiency of the algorithm can be further elaborated by an estimate of the number of iterations for the guaranteed convergence.  The convergence condition of NST+HT+FB is also improved.  Moreover, an adaptive scheme (AdptNST+HT+FB) without the knowledge of the sparsity level is proposed with its convergence guarantee.  The number of iterations for the finite step of convergence of the AdptNST+HT+FB scheme is also derived. It is further shown that the number of iterations can be significantly reduced by exploiting the structure of the specific sparse signal or the random measurement matrix.
\end{abstract}

\MOS{<Subject classification numbers>}


\maketitle


\vspace{-6pt}

\section{Introduction}
Compressive sensing (CS) is one of the most relevant techniques of signal sampling and reconstruction in modern data sciences.  The main aim is to recover sparse signals from incomplete linear measurements
\begin{equation}\label{yisAx}
y=Ax,
\end{equation}
where $A\in\mathbb{R}^{M\times N}$ is the sampling matrix with $M\ll N$, and $x$ denotes the $N$-dimensional sparse signal with only $s$ nonzero
coefficients.

Since most natural signals are sparse or highly compressible under a basis, CS has a wide range of applications including signal processing \cite{R1}, sensor network \cite{R2}, biological application \cite{R3}, sub-Nyquist sampling system \cite{R4}, etc. Various algorithms have been proposed for solving problem (\ref{yisAx}). Evidently, the underlying model involves finding the sparsest solutions satisfying the linear equations,
\begin{equation}\label{l0formulation}
\min \limits_{x\in\mathbb{R}^{N}}\|x\|_{0},\ \ s.t.\ \ y=Ax,
\end{equation}
or, one of its Lagrangian versions
\begin{equation}\label{l0Lagrange}
\min \limits_{x\in\mathbb{R}^{N}}\|y-Ax\|_{2}^{2}+\lambda\|x\|_{0},
\end{equation}
where $\|x\|_{0}$ is $\ell_{0}$ ``norm'' of the vector $x\in\mathbb{R}^{N}$, indicating the number of nonzero entries in $x$ and $\lambda$ is a regularization parameter. (\ref{l0formulation}) and (\ref{l0Lagrange}) are clearly combinatorial and computationally intractable \cite{R5}.

Heuristic (tractable) algorithms have been extensively studied to solve this problem. A common strategy is based on problem relaxation that replaces the $\ell_{0}$ ``norm'' by an $\ell_{p}$-norm (with $0<p\leq1$). Relaxed problems can be solved efficiently by simple optimization procedures. Well-known algorithms using the strategy are basis pursuit \cite{R6,R7,R8}, $\ell_{p}$ optimization model \cite{R9,R10,R11}, iterative reweighted $\ell_{1}$-norm \cite{R12}, least absolute shrinkage and selection operator \cite{R13}, focal underdetermined system solver \cite{R14}, bregman iterations \cite{R15}. The pursuit algorithm is another popular category, which builds up the sparse solutions by making a series of greedy decisions. Typical representative approaches are matching pursuit \cite{R16}, orthogonal matching pursuit \cite{R17}. In addition, a number of variants of the greedy pursuit algorithms have also been proposed, e.g., stagewise orthogonal matching pursuit \cite{R18}, compressive sampling matching pursuit \cite{R19} and subspace pursuit \cite{R20}, etc.

In addition, a common approach is the iterative thresholding/shrinkage algorithm, which has attracted tremendous amount of attention due to its remarkable performance/complexity trade-offs. These methods recover the sparse signal by making a succession of thresholding operations. The iterative hard thresholding algorithm was first introduced in \cite{R21}. Following a similar procedure, a soft thresholding mechanism was proposed in \cite{R22}.
\newpage\noindent  A hard thresholding pursuit algorithm was also a popular procedure \cite{R23,R24}. In fact, the hard thresholding pursuit algorithm can be regarded as a hybrid of  the iterative hard thresholding algorithm and the compressive sampling matching pursuit. In \cite{R25}, a thresholding with feedback and null-space turning (NST+HT+FB) algorithm was proposed to find sparse solutions. The proposed algorithms were brought into a concise framework of null space tuning (NST). Several sparsity enhancing operators were incorporated into the NST framework to develop various algorithms. These algorithms were shown to be exceedingly fast and effective, particularly for large scale systems.

As shown in the present article, the NST+HT+FB algorithm converges to the true solution under a certain  (preconditioned) restricted isometry condition. In this paper, it turns out that further efforts at understanding the NST+HT+FB algorithm reveal that the finite number of convergence steps can be explicitly estimated.  An adaptive NST+HT+FB (AdptNST+HT+FB) procedure without the knowledge of the sparsity level is further investigated.  Analysis of the finite convergence of the AdptNST+HT+FB is also carried out.

For clarity, notations are used as follows in this article. For any $c\in\mathbb{R}$, $\lceil c\rceil$ is the smallest integer that is greater than or equal to $c$. $S$ is the support of $s$-sparse vector $x$. $x_{T}$ is the restriction of a vector $x\in\mathbb{R}^{N}$ to an index set $T$. We denote by $T^{c}$ the complement set of $T$ in $\{1,2,\ldots,N\}$, and by $A_{T}$ the sub-matrix consisting of columns of $A$ indexed
by $T$, respectively. $A_{l}$ denotes the $l$th column of the matrix $A$. $\widehat{x}\in\mathbb{R}^{N}_{+}$ is the nonincreasing rearrangement of a vector $x=(x_1, x_2,\cdots,x_N)^{'}\in\mathbb{R}^{N}$, i.e., $\widehat{x}_{1}\geq \widehat{x}_{2}\geq\ldots \widehat{x}_{N}\geq0$
and there exists a permutation $\pi$ of $\{1,\ldots,N\}$ such that $\widehat{x}_{i}=|x_{\pi(i)}|$ for all $i\in\{1,\ldots,N\}$.
$T\triangle T'$ is the symmetric difference of $T$ and $T'$, i.e., $T\triangle T' =(T\setminus T')\cup(T'\setminus T)$. $|T|$ is the cardinality of set $T$.
\section{The framework of NST+HT+FB}
The iterative framework of the approximation and null space tuning (NST) algorithms is as follows
\begin{equation*} \label{eq:1}
\left\{ \begin{aligned}
         \begin{aligned}
         &u^{k}=\mathbb{D}(x^{k}),\\
         &x^{k+1}=x^{k}+\mathbb{P}(u^{k}-x^{k}),\\
         \end{aligned}
         \end{aligned} \right.
\end{equation*}
where $\mathbb{D}(x^{k})$ approximates the desired sparse solution by various principles (Here $\mathbb{D}$ is set as thresholding plus a feedback), and $\mathbb{P}:=I-A^{\ast}(AA^{\ast})^{-1}A$ is the orthogonal projection
onto ker$A$. The feasibility of $x^{0}$ is assumed, which guarantees that the sequence $\{x^{k}\}$ are all feasible. Obviously, $u^{k}\rightarrow x$ is
expected as $k$ increases. Since the sequence $\{x^{k}\}$ are always feasible in the framework of the NST algorithms, one may split $y$ as
\begin{equation*}
y=Ax=A_{T_{k}}x_{T_{k}}^{k}+A_{T^{c}_{k}}x_{T^{c}_{k}}^{k}.
\end{equation*}
In most (if not all) thresholding algorithms, thresholding (hard or soft) is taken by merely keeping the entries of $x^{k}$ on $T_{k}$,
and completely abandons the contribution of $A_{T^{c}_{k}}x_{T^{c}_{k}}^{k}$ to the measurement $y$. Though $x_{T^{c}_{k}}^{k}$ gradually diminishes as $k\rightarrow\infty$, it is not difficult to observe that the contribution of $A_{T^{c}_{k}}x_{T^{c}_{k}}^{k}$ to $y$ can be quite significant at initial
iterations. Therefore, simple thresholding alone can be quite infeasible at earlier stages. The mechanism of feedback is to feed the contribution of $A_{T_{k}^{c}}x_{T_{k}^{c}}^{k}$ to $y$ back to im($A_{T_{k}}$), the image of $A_{T_{k}}$. One straightforward way is to set
\begin{equation}\label{eq3}
\eta^{k}=\arg\min \limits_{\eta}\|A_{T_{k}}\eta-A_{T_{k}^{c}}x_{T_{k}^{c}}^{k}\|_{2},
\end{equation}
which has the best/least-square solution
\begin{equation*}
\eta^{k}=(A_{T_{k}}^{\ast}A_{T_{k}})^{-1}A_{T_{k}}^{\ast}A_{T_{k}^{c}}x_{T_{k}^{c}}^{k}.
\end{equation*}
The NST+HT+FB algorithm is then established as follows
\begin{equation*} \label{eq:1}
\text{(NST+HT+FB)}\ \ \ \
\left\{ \begin{aligned}
         \begin{aligned}
         &\mu_{T_{k}}^{k}=x_{T_{k}}^{k}+(A_{T_{k}}^{\ast}A_{T_{k}})^{-1}A_{T_{k}}^{\ast}A_{T_{k}^{c}}x_{T_{k}^{c}}^{k}, \\
         &\mu_{T_{k}^{c}}^{k}=0, \\
         &x^{k+1}=x^{k}+\mathbb{P}(u^{k}-x^{k}).\\
         \end{aligned}
         \end{aligned} \right.
\end{equation*}
Since $|T_{k}|=s$ at each iteration, NST+HT+FB constructs a sequence $\{\mu^{k}\}$ of $s$-sparse signals. With $\mathbb{P}:=I-A^{\ast}(AA^{\ast})^{-1}A$, the null space tuning (NST) step can be rewritten as $x^{k+1}=u^{k}+A^{\ast}(AA^{\ast})^{-1}(y-Au^{k})$.

 Let $x$ be the solution to $y=Ax$ with only $s$ sparsity. \cite{R25} shows that if the preconditioned restricted isometry property (P-RIP) and restricted isometry property (RIP) constants of $A$ satisfies $\delta_{2s}+\sqrt{2}\gamma_{3s}<1$, then $u^{k}$ in NST+HT+FB converges to $x$. This paper presents an improved convergence condition $\delta_{2s}^{2}+2\gamma_{3s}^{2}<1$ and the number of the finite step for the convergence of NST+HT+FB is also estimated.

 The standard NST+HT+FB algorithm requires the knowledge of the sparsity level of the desired solution.  This seems wishful in most applications.  Another adaptive scheme of NST+HT+FB (AdptNST+HT+FB) is introduced in this work.   The AdptNST+HT+FB avoids the prior estimation of the sparsity level, in which the sparsity level is adjusted upward gradually. Specifically, a sequence $\{\mu^{k}\}$ of $k$-sparse signals is established according to
 \begin{equation*} \label{eq:1}
\text{(AdptNST+HT+FB)}\ \ \ \
\left\{ \begin{aligned}
         \begin{aligned}
         &\mu_{T_{k}}^{k}=x_{T_{k}}^{k}+(A_{T_{k}}^{\ast}A_{T_{k}})^{-1}A_{T_{k}}^{\ast}A_{T_{k}^{c}}x_{T_{k}^{c}}^{k}, \\
         &\mu_{T_{k}^{c}}^{k}=0, \\
         &x^{k+1}=x^{k}+\mathbb{P}(u^{k}-x^{k}),\\
         \end{aligned}
         \end{aligned} \right.
\end{equation*}
where $|T_{k}|=k$ at each iteration. The convergence condition and the number of iterations for the AdptNST+HT+FB scheme are also fully studied and presented.
\section{Preliminary results}

\begin{definition}\label{RIP} \cite{R26}. For each integer $s = 1, 2,... ,$ the restricted isometry constant $\delta_{s}$ of a matrix $A$ is defined as the smallest number $\delta_{s}$ such that
\begin{equation*}
(1-\delta_{s})\|x\|_{2}^{2}\leq\|Ax\|_{2}^{2}\leq(1+\delta_{s})\|x\|_{2}^{2}
\end{equation*}
holds for all $s$-sparse vectors $x$. Equivalently, it is given by
\begin{equation*}
\delta_{s}=\max_{|S|\leq s}\|I-A_{S}^{\ast}A_{S}\|_{2\rightarrow2}.
\end{equation*}
\end{definition}
\begin{definition}\label{P-RIP} \cite{R25}. For each integer $s = 1, 2,... ,$ the preconditioned restricted isometry constant $\gamma_{s}$ of a matrix $A$ is defined as the smallest number $\gamma_{s}$ such that
\begin{equation*}
(1-\gamma_{s})\|x\|_{2}^{2}\leq\|(AA^{\ast})^{-\frac{1}{2}}Ax\|_{2}^{2}
\end{equation*}
holds for all $s$-sparse vectors $x$.
In fact, the preconditioned restricted isometry constant $\gamma_{s}$ characterizes the restricted isometry property of the preconditioned matrix $(AA^{\ast})^{-\frac{1}{2}}A$. Since
\begin{equation*}
\|(AA^{\ast})^{-\frac{1}{2}}Ax\|_{2}\leq\|(AA^{\ast})^{-\frac{1}{2}}A\|_{2}\|x\|_{2}=\|x\|_{2},
\end{equation*}
$\gamma_{s}$ is actually the smallest number such that, for all $s$-sparse vectors $x$,
\begin{equation*}
(1-\gamma_{s})\|x\|_{2}^{2}\leq\|(AA^{\ast})^{-\frac{1}{2}}Ax\|_{2}^{2}\leq(1+\gamma_{s})\|x\|_{2}^{2}.
\end{equation*}
It indicates $\gamma_{s}(A)=\delta_{s}((AA^{\ast})^{-\frac{1}{2}}A)$. Equivalently, it is given by
\begin{equation*}
\gamma_{s}=\max_{|S|\leq s}\|I-A_{S}^{\ast}(AA^{\ast})^{-1}A_{S}\|_{2\rightarrow2}.
\end{equation*}
\end{definition}

\begin{lemma} \label{ert}
Let $\delta_t$ be the RIP constant of $A$. For $u,v\in\mathbb{R}^{N}$, if $|\supp (u) \cup \supp (v)|\leq t$, then $|\langle u,(I-A^{\ast}A)v\rangle|\leq\delta_{t}\|u\|_{2}\|v\|_{2}$. Suppose $|T' \cup \supp (v)|\leq t$, then $\|((I-A^{\ast}A)v)_{T'}\|_{2}\leq\delta_{t}\|v\|_{2}$.
\end{lemma}

$\bm{Proof.}$
Setting $T=\supp (v)\cup  \supp (u)$, $|T|\leq t$, one has
\begin{equation*}
\begin{aligned}
|\langle u,(I-A^{\ast}A)v\rangle| & =|\langle u,v\rangle-\langle Au,Av\rangle| \\
&=|\langle u_{T},v_{T}\rangle-\langle A_{T}u_{T},A_{T}v_{T}\rangle|\\
&=|\langle u_{T},(I-A_{T}^{\ast}A_{T})v_{T}\rangle|\\
&\leq\|u_{T}\|_{2}\|(I-A_{T}^{\ast}A_{T})v_{T}\|_{2}\\
&\leq\|u_{T}\|_{2}\|I-A_{T}^{\ast}A_{T}\|_{2\rightarrow2}\|v_{T}\|_{2}\\
&\leq\delta_{t}\|u\|_{2}\|v\|_{2}.\\
\end{aligned}
\end{equation*}
The first inequality is due to the Cauchy-Schwarz inequality and the second inequality is due to the submultiplicativity of matrix norms, while the last step is based on {\em Definition \ref{RIP}}.
Since
\begin{equation*}
\begin{aligned}
\|((I-A^{\ast}A)v)_{T'}\|_{2}^{2} &=\langle((I-A^{\ast}A)v)_{T'},(I-A^{\ast}A)v\rangle\\
&\leq\delta_{t}\|((I-A^{\ast}A)v)_{T'}\|_{2}\|v\|_{2},
\end{aligned}
\end{equation*}
one can obtain that $\|((I-A^{\ast}A)v)_{T'}\|_{2}\leq\delta_{t}\|v\|_{2}$.

\begin{remark} \label{nin}
Let $\gamma_{t}$ be the P-RIP constant of $A$, i.e., $\gamma_{t}(A)=\delta_{t}((AA^{\ast})^{-\frac{1}{2}}A)$. For $u,v\in\mathbb{R}^{N}$, if $|\supp (u) \cup \supp (v)|\leq t$, then $|\langle u,(I-A^{\ast}(AA^{\ast})^{-1}A)v\rangle|\leq\gamma_{t}\|u\|_{2}\|v\|_{2}$. Suppose $|T' \cup \supp (v)|\leq t$, then $\|((I-A^{\ast}(AA^{\ast})^{-1}A)v)_{T'}\|_{2}\leq\gamma_{t}\|v\|_{2}$.
\end{remark}
%

\begin{lemma}\label{lemm}
For $e\in\mathbb{R}^{M}$, $\|(A^{\ast}e)_{T}\|_{2}\leq\sqrt{1+\delta_{t}}\|e\|_{2}$, when $|T|\leq t$.
\end{lemma}

$\bm{Proof.}$
\begin{equation*}
\begin{aligned}
\|(A^{\ast}e)_{T}\|_{2}^{2}&=\langle A^{\ast}e,(A^{\ast}e)_{T}\rangle=\langle e,A(A^{\ast}e)_{T}\rangle\\
&\leq\|e\|_{2}\|A(A^{\ast}e)_{T}\|_{2}\leq\|e\|_{2}\sqrt{1+\delta_{t}}\|(A^{\ast}e)_{T}\|_{2},
\end{aligned}
\end{equation*}
Hence, for all $e\in\mathbb{R}^{M}$, we have $\|(A^{\ast}e)_{T}\|_{2}\leq\sqrt{1+\delta_{t}}\|e\|_{2}$.

\begin{remark}\label{rrrrrrr}
If $\delta_{s}((AA^{\ast})^{-1}A)=\theta_{s}$, then for $e\in\mathbb{R}^{M}$, $\|(A^{\ast}(AA^{\ast})^{-1}e)_{T}\|_{2}\leq\sqrt{1+\theta_{t}}\|e\|_{2}$, when $|T|\leq t$.
\end{remark}

\section{Related preliminaries}
In this section, we introduce basic preliminaries that will be used later. They demonstrate the closeness of $\mu^{k}$ to $x$ and state the fact that how the indices of nonzero entries of $x$ are captured in the sequences produced by NST+HT+FB and AdptNST+HT+FB.
\begin{lemma}\label{h}
 Suppose $y=Ax+e$ where $x\in\mathbb{R}^{N}$ is $s$-sparse with $S=$\supp$(x)$ and $e\in\mathbb{R}^{M}$ is the measurement error. If $\mu'\in\mathbb{R}^{N}$ is $s'$-sparse and $T$ is an index set of $t\geq s$ largest absolute entries of $\mu'+A^{\ast}(AA^\ast)^{-1}(y-A\mu')$, then
\begin{equation*}
\|x_{T^{c}}\|_{2}\leq\sqrt{2}(\gamma_{s+s'+t}\|x-\mu'\|_{2}+\sqrt{1+\theta_{t+s}}\|e\|_{2}),
\end{equation*}
where $\theta_{s}(A)=\delta_{s}((AA^\ast)^{-1}A)$.
\end{lemma}

$\bm{Proof.}$
We first have that

$\|[\mu'+A^{\ast}(AA^\ast)^{-1}(y-A\mu')]_{T}\|_{2}\geq\|[\mu'+A^{\ast}(AA^\ast)^{-1}(y-A\mu')]_{S}\|_{2}.$

Eliminating the common terms over $T\bigcap S$, one has

$\|[\mu'+A^{\ast}(AA^\ast)^{-1}(y-A\mu')]_{T\setminus S}\|_{2}\geq\|[\mu'+A^{\ast}(AA^\ast)^{-1}(y-A\mu')]_{S\setminus T}\|_{2}.$

For the left hand,
\begin{equation*}
\begin{aligned}
&\|[\mu'+A^{\ast}(AA^\ast)^{-1}(y-A\mu')]_{T\setminus S}\|_{2}\\
&=\|[\mu'-x+A^{\ast}(AA^\ast)^{-1}(Ax+e-A\mu')]_{T\setminus S}\|_{2}\\
&=\|[(I-A^{\ast}(AA^\ast)^{-1}A)(\mu'-x)+A^{\ast}(AA^\ast)^{-1}e]_{T\setminus S}\|_{2}.\\
\end{aligned}
\end{equation*}

The right hand satisfies
\begin{equation*}
\begin{aligned}
&\|[\mu'+A^{\ast}(AA^\ast)^{-1}(y-A\mu')]_{S\setminus T}\|_{2}\\
&=\|[\mu'+A^{\ast}(AA^\ast)^{-1}(Ax+e-A\mu')+x-x]_{S\setminus T}\|_{2}\\
&\geq\|x_{S\setminus T}\|_{2}-\|[(I-A^{\ast}(AA^\ast)^{-1}A)(\mu'-x)+A^{\ast}(AA^\ast)^{-1}e]_{S\setminus T}\|_{2}.\\
\end{aligned}
\end{equation*}
Consequently,
\begin{equation*}
\begin{aligned}
&\|x_{S\setminus T}\|_{2}\leq \|[(I-A^{\ast}(AA^\ast)^{-1}A)(\mu'-x)+A^{\ast}(AA^\ast)^{-1}e]_{S\setminus T}\|_{2}\\
&+\|[(I-A^{\ast}(AA^\ast)^{-1}A)(\mu'-x)+A^{\ast}(AA^\ast)^{-1}e]_{T\setminus S}\|_{2}\\
&\leq\sqrt{2}\|[(I-A^{\ast}(AA^\ast)^{-1}A)(\mu'-x)+A^{\ast}(AA^\ast)^{-1}e]_{T\triangle S}\|_{2}\\
&\leq\sqrt{2}\|[(I-A^{\ast}(AA^\ast)^{-1}A)(\mu'-x)]_{T\triangle S}\|_{2}+\sqrt{2}\|[A^{\ast}(AA^\ast)^{-1}e]_{T\triangle S}\|_{2}\\
&\leq\sqrt{2}(\gamma_{s+s'+t}\|x-\mu'\|_{2}+\sqrt{1+\theta_{t+s}}\|e\|_{2}).
\end{aligned}
\end{equation*}
 The last step is due to {\em Remark \ref{nin}} and {\em Remark \ref{rrrrrrr}}.

\begin{lemma}\label{lemm6}
 Suppose $y=Ax+e$ where $x\in\mathbb{R}^{N}$ is $s$-sparse with $S=$\supp$(x)$ and $e\in\mathbb{R}^{M}$ is the measurement error. Let $T=$\supp$(x')$ and $|T|=t$. If $\mu'$ is the feedback of $x'$ that subjects to $\mu_{T}'=x_{T}'+(A_{T}^{\ast}A_{T})^{-1}A_{T}^{\ast}A_{T^{c}}x_{T^{c}}'$ and $\mu_{T^{c}}'=0$, then
\begin{equation*}
\|(x-\mu')\|_{2}\leq\frac{\sqrt{1+\delta_{t}}\|e\|_{2}}{1-\delta_{s+t}}+\frac{\|x_{T^{c}}\|_{2}}{\sqrt{1-\delta_{s+t}^{2}}}.
\end{equation*}
\end{lemma}

$\bm{Proof.}$
For any $z\in\mathbb{R}^{N}$ supported on $T$,
\begin{equation*}
\begin{aligned}
&\langle A\mu'-y,Az\rangle\\
&=\langle A_{T}x_{T}'+A_{T}( A_{T}^{\ast} A_{T})^{-1}A_{T}^{\ast} A_{T^{c}}x_{T^{c}}'-y,A_{T}z_{T}\rangle\\
&=\langle A_{T}^{\ast}(A_{T}x_{T}'+A_{T^{c}}x_{T^{c}}'-y),z_{T}\rangle\\
&=\langle A_{T}^{\ast}(Ax'-y),z_{T}\rangle\\
&=0\\
\end{aligned}
\end{equation*}
The last step is due to the feasibility of $x'$, i.e., $y=Ax'$.
The inner product can also be written as

$\langle A\mu'-y,Az\rangle=\langle (A\mu'-Ax-e),Az\rangle=0$.

\noindent Therefore,

$\langle (\mu'-x),A^{\ast}Az\rangle=\langle e,Az\rangle,~~\forall z\in\mathbb{R}^{N}$ supported on $T$.

\noindent Since $(\mu'-x)_{T}$ is supported on $T$, one has

$\langle (\mu'-x),A^{\ast}A(\mu'-x)_{T}\rangle=\langle e,A(\mu'-x)_{T}\rangle.$

\noindent Consequently,
\begin{equation*}
\begin{aligned}
&\|(\mu'-x)_{T}\|_{2}^{2}=\\
&=\langle(\mu'-x),(\mu'-x)_{T}\rangle\\
&=|\langle (x-\mu'),(I-A^{\ast}A)(x-\mu')_{T}\rangle+\langle e,A(\mu'-x)_{T}\rangle|\\
&\leq\delta_{s+t}\|x-\mu'\|_{2}\|(x-\mu')_{T}\|_{2}+\sqrt{1+\delta_{t}}\|e\|_{2}\|(x-\mu')_{T}\|_{2}\\
\end{aligned}
\end{equation*}
Using  {\em Lemma \ref{ert}}, Cauchy-Schwarz inequality and {\em Definition \ref{RIP}} can obtain the last inequality.
\noindent Therefore, we have

$\|(x-\mu')_{T}\|_{2}\leq\delta_{s+t}\|x-\mu'\|_{2}+\sqrt{1+\delta_{t}}\|e\|_{2}$.

\noindent It then follows that

\begin{equation*}
\begin{aligned}
&\|(x-\mu')\|_{2}^{2}=\|(x-\mu')_{T}\|_{2}^{2}+\|(x-\mu')_{T^{c}}\|_{2}^{2}\\
&\leq(\delta_{s+t}\|x-\mu'\|_{2}+\sqrt{1+\delta_{t}}\|e\|_{2})^{2}+\|x_{T^{c}}\|_{2}^{2}.\\
\end{aligned}
\end{equation*}

\noindent In other words,

$(\sqrt{1-\delta_{s+t}^{2}}\|(x-\mu')\|_{2}-\frac{\delta_{s+t}\sqrt{1+\delta_{t}}}{\sqrt{1-\delta_{s+t}^{2}}}\|e\|_{2})^{2}
\leq\frac{1+\delta_{t}}{1-\delta_{s+t}^{2}}\|e\|_{2}^{2}+\|x_{T^{c}}\|_{2}^{2}.$

\noindent It means that
\begin{equation*}
\begin{aligned}
&\|(x-\mu')\|_{2}\\
&\leq\frac{\delta_{s+t}\sqrt{1+\delta_{t}}\|e\|_{2}+\sqrt{(1+\delta_{t})\|e\|_{2}^{2}+({1-\delta_{s+t}^{2})\|x_{T^{c}}\|_{2}^{2}}}}{1-\delta_{s+t}^{2}}\\
&\leq\frac{\|x_{T^{c}}\|_{2}}{\sqrt{1-\delta_{s+t}^{2}}}+\frac{\sqrt{1+\delta_{t}}\|e\|_{2}}{1-\delta_{s+t}}.\\
\end{aligned}
\end{equation*}

\begin{corollary}\label{ha}
 Let $x\in\mathbb{R}^{N}$ be $s$-sparse and $y=Ax+e$ for $e\in\mathbb{R}^{M}$. If $\{u^{k}\}$ is the sequence of AdptNST+HT+FB, then $u^{k}$ satisfies
\begin{equation*}\label{hhh}
\begin{aligned}
&\|(x-u^{k})\|_{2}\\
&\leq\sqrt{\frac{2\gamma_{s+2k-1}^{2}}{(1-\delta_{s+k}^{2})}}\|x-\mu^{k-1}\|_{2}\\
&+(\frac{\sqrt{1+\delta_{k}}}{1-\delta_{s+k}}+
\frac{\sqrt{2(1+\theta_{s+k})}}{\sqrt{1-\delta_{s+k}^{2}}})\|e\|_{2},~k\geq s.\\
\end{aligned}
\end{equation*}
\end{corollary}

$\bm{Proof.}$
 Applying {\em Lemma \ref{h}} to $\mu'=\mu^{k-1}$ and $T=T_{k}$ for $k\geq s$ gives
\begin{equation*}
\|x_{T^{c}_{k}}\|\leq\sqrt{2}(\gamma_{s+2k-1}\|x-\mu^{k-1}\|_{2}+\sqrt{1+\theta_{s+k}}\|e\|_{2}),
\end{equation*}
\noindent and setting $\mu'=u^{k}$ and $T=T_{k}$ in {\em Lemma \ref{lemm6}} obtains
\begin{equation*}
\|(x-u^{k})\|_{2}\leq\frac{\|x_{T^{c}_{k}}\|_{2}}{\sqrt{(1-\delta_{s+k}^{2})}}+\frac{\sqrt{1+\delta_{k}}\|e\|_{2}}{1-\delta_{s+k}}.
\end{equation*}
\noindent Combining these two inequalities, we have
\begin{equation*}
\begin{aligned}
&\|(x-u^{k})\|_{2}\\
&\leq\sqrt{\frac{2\gamma_{s+2k-1}^{2}}{(1-\delta_{s+k}^{2})}}\|x-\mu^{k-1}\|_{2}\\
&+(\frac{\sqrt{1+\delta_{k}}}{1-\delta_{s+k}}+
\frac{\sqrt{2(1+\theta_{s+k})}}{\sqrt{1-\delta_{s+k}^{2}}})\|e\|_{2}.
\end{aligned}
\end{equation*}

Through {\em Corollary \ref{ha}}, if the P-RIP and RIP constants of $A$ satisfies $2\gamma_{s+2k-1}^{2}+\delta_{s+k}^{2}<1$, then the sequence of $\{\mu_{k}\}$ in AdptNST+HT+FB converges to $x$. If the prior estimation of the sparsity is known $(|T_{k}|=s)$, we can have the following remark for NST+HT+FB.
\begin{remark}\label{hhhh1}
If $\{u^{k}\}$ is the sequence of NST+HT+FB, then $u^{k}$ satisfies
\begin{equation*}
\begin{aligned}
\|x-u^{k}\|_{2}&\leq\sqrt{\frac{2\gamma_{3s}^{2}}{(1-\delta_{2s}^{2})}}\|x-\mu^{k-1}\|_{2}\\
&+(\frac{\sqrt{1+\delta_{s}}}{1-\delta_{2s}}+
\frac{\sqrt{2(1+\theta_{2s})}}{\sqrt{1-\delta_{2s}^{2}}})\|e\|_{2},~k\geq1.
\end{aligned}
\end{equation*}
\end{remark}

As shown in {\em Remark \ref{hhhh1}}, if the P-RIP and RIP constants of $A$ satisfies $\delta_{2s}^{2}+2\gamma_{3s}^{2}<1$, then the sequence of $\{u^{k}\}$ in NST+HT+FB converges to $x$. Compared to the condition $\delta_{2s}+\sqrt{2}\gamma_{3s}<1$ in \cite{R25}, the condition in {\em Remark \ref{hhhh1}} is obvious improved. Furthermore, if $A$ is the Parseval frame, the P-RIP and RIP condition is relaxed to RIP condition, i.e., $\delta_{2s}^{2}+2\delta_{3s}^{2}<1$.
\section{The number of iterations for convergence of NST+HT+FB}
This section contains the main result about NST+HT+FB that how many iterations are necessary to correctly capture the support of $x$. Besides the general sparse signal recovery that the measurement matrix can be used to all sparse signals simultaneously, some nonuniform cases about the necessary iterations by exploiting the extra information of the sparse signal are also shown. For notational simplicity, we define  $\rho_{s}=\sqrt{\frac{2\gamma_{s}^{2}}{1-\delta_{s}^{2}}}$, $\tau_{s}=(\frac{\sqrt{1+\delta_{s}}}{1-\delta_{s}}+
\frac{\sqrt{2(1+\theta_{s})}}{\sqrt{1-\delta_{s}^{2}}})$, and $\omega_{s}=\frac{\sqrt{2}\gamma_{s}\rho_{s}^{l-1}\sqrt{1+\delta_{s}}}{1-\delta_{s}}
+\frac{\sqrt{2}\gamma_{s}\tau_{s}(1-\rho_{s}^{l-1})}{1-\rho_{s}}+\sqrt{2(1+\theta_{s})}$.

\subsection{Uniform sparse recovery}
Due the mechanism of feedback, the algorithm converges when $T_{k}=S$. The remaining topic is to find the number of steps needed for capturing the true support $S$. The following lemmas show the size of the indices of nonzero entries of $x$ that captured in the support sets $\{T_{k}\}$ produced by NST+HT+FB increases by the iteration gradually. The number of iterations for increasing a specified amount of correct indices are derived.
\begin{lemma}\label{qw}
Suppose $y=Ax+e$ with an $s$-sparse vector $x$ and let $\{T_{k}\}$ be index sets of the sequence $\{u^{k}\}$ in NST+HT+FB. For integers $k,p\geq0$, if $T_{k}$ contains the indices of $p$ largest absolute entries of $x$. Then $\{T_{k+l}\}$ contains the indices of $p+q$ largest absolute entries of $x$ via $l$ iterations for integers $l,q\geq1$, provided
\begin{equation*}
\widehat{x}_{p+q}>\rho_{3s}^{l}\|\widehat{x}_{\{p+1,\ldots,s\}}\|_{2}+\omega_{3s}\|e\|_{2}.
\end{equation*}
\end{lemma}

$\bm{Proof.}$
For NST+HT+FB, the hypothesis is $\pi(\{1,\ldots,p\})\subseteq T_{k}$ and the goal is to prove that $\pi(\{1,\ldots,p+q\})\subseteq T_{k+l}$. That is to say the $|(\mu^{l+k-1}+A^{\ast}(AA^{\ast})^{-1}(y-A\mu^{l+k-1}))_{\pi(j)}|$ for $j\in\{1,\ldots p+q\}$ are among the $s$ largest entries of
$|(\mu^{l+k-1}+A^{\ast}(AA^{\ast})^{-1}(y-A\mu^{l+k-1}))_{i}|$ for $i\in\{1,\ldots,N\}$. Since \supp$(x)=S$ and $|S|$=s, it is the enough to prove that
\begin{equation*}
\begin{aligned}
&\min \limits_{j\in\{1,\ldots,p+q\}} |(\mu^{l+k-1}+A^{\ast}(AA^{\ast})^{-1}(y-A\mu^{l+k-1}))_{\pi(j)}|\\
&>\max \limits_{i\in S^{c}} |(\mu^{l+k-1}+A^{\ast}(AA^{\ast})^{-1}(y-A\mu^{l+k-1}))_{i}|.\\
\end{aligned}
\end{equation*}
\noindent For $j\in\{1,\ldots p+q\}$ and $i\in S^{c}$ , in view of
\begin{equation*}
\begin{aligned}
&|(\mu^{l+k-1}+A^{\ast}(AA^{\ast})^{-1}(y-A\mu^{l+k-1}))_{\pi(j)}|\\
&\geq|x_{\pi(j)}|-|(-x+\mu^{l+k-1}+A^{\ast}(AA^{\ast})^{-1}(y-A\mu^{l+k-1}))_{\pi(j)}|\\
&\geq \widehat{x}_{p+q}-|((I-A^{\ast}(AA^{\ast})^{-1}A)(\mu^{l+k-1}-x)+A^{\ast}(AA^{\ast})^{-1}e)_{\pi(j)}|,\\
&|(\mu^{l+k-1}+A^{\ast}(AA^{\ast})^{-1}(y-A\mu^{l+k-1}))_{i}|\\
&=|(-x+\mu^{l+k-1}+A^{\ast}(AA^{\ast})^{-1}(y-A\mu^{l+k-1}))_{i}|\\
&=|((I-A^{\ast}(AA^{\ast})^{-1}A)(\mu^{l+k-1}-x)+A^{\ast}(AA^{\ast})^{-1}e)_{i}|,
\end{aligned}
\end{equation*}
\noindent one only needs to prove next, for all $j\in\{1,\ldots,p+q\}$ and $i\in S^{c}$,
\begin{equation*}
\begin{aligned}
\widehat{x}_{p+q}&>|[(I-A^{\ast}(AA^{\ast})^{-1}A)(\mu^{l+k-1}-x)+A^{\ast}(AA^{\ast})^{-1}e]_{\pi(j)}|\\
&+|((I-A^{\ast}(AA^{\ast})^{-1}A)(\mu^{l+k-1}-x)+A^{\ast}(AA^{\ast})^{-1}e)_{i}|.\\
\end{aligned}
\end{equation*}
\noindent The right hand side can be bounded by
\begin{equation*}
\begin{aligned}
&\sqrt{2}\|[(I-A^{\ast}(AA^{\ast})^{-1}A)(\mu^{l+k-1}-x)+A^{\ast}(AA^{\ast})^{-1}e]_{\{\pi(j),i\}}\|_{2}\\
&\leq\sqrt{2}\|[(I-A^{\ast}(AA^{\ast})^{-1}A)(\mu^{l+k-1}-x)]_{\{\pi(j),i\}}\|_{2}\\
&+\sqrt{2}\|[A^{\ast}(AA^{\ast})^{-1}e]_{\{\pi(j),i\}}\|_{2}\\
&\leq\sqrt{2}\gamma_{2s+2}\|\mu^{l+k-1}-x\|_{2}+\sqrt{2(1+\theta_{2})}\|e\|_{2}\\
&\leq\sqrt{2}\gamma_{2s+2}(\rho_{3s}^{l-1}\|u^{k}-x\|_{2}+\frac{\tau_{2s}(1-\rho_{3s}^{l-1})\|e\|_{2}}{1-\rho_{3s}})\\
&+\sqrt{2(1+\theta_{2})}\|e\|_{2},\\
\end{aligned}
\end{equation*}
where {\em Remark \ref{hhhh1}} was used $l-1$ times in the last step. From {\em Lemma \ref{lemm6}} and the assumption that $\pi(\{1,\ldots,p\})\subseteq T_{k}$,
\begin{equation*}
\begin{aligned}
&\sqrt{2}\|[(I-A^{\ast}(AA^{\ast})^{-1}A)(\mu^{l+k-1}-x)+A^{\ast}(AA^{\ast})^{-1}e]_{\{\pi(j),i\}}\|_{2}\\
&\leq\sqrt{2}\gamma_{2s+2}\rho_{3s}^{l-1}(\frac{1
}{\sqrt{(1-\delta_{2s}^{2})}}\|x_{T^{c}_{k}}\|_{2}+\frac{\sqrt{1+\delta_{s}}}{1-\delta_{2s}}\|e\|_{2})\\
&+\frac{\sqrt{2}\gamma_{2s+2}\tau_{2s}(1-\rho_{3s}^{l-1})\|e\|_{2}}{1-\rho_{3s}}+\sqrt{2(1+\theta_{2})}\|e\|_{2}\\
&\leq\rho_{3s}^{l}\|x_{T^{c}_{k}}\|_{2}\\
&+(\sqrt{2}\gamma_{2s+2}\rho_{3s}^{l-1}\frac{\sqrt{1+\delta_{s}}}{1-\delta_{2s}}+\frac{\sqrt{2}\gamma_{2s+2}\tau_{2s}(1-\rho_{3s}^{l-1})}{1-\rho_{3s}}\\
&+\sqrt{2(1+\theta_{2})})\|e\|_{2}\\
&\leq\rho_{3s}^{l}\|x_{\pi(\{1,\ldots,p\})^{c}}\|_{2}\\
&+(\sqrt{2}\gamma_{3s}\rho_{3s}^{l-1}\frac{\sqrt{1+\delta_{3s}}}{1-\delta_{3s}}+\frac{\sqrt{2}\gamma_{3s}\tau_{3s}(1-\rho_{3s}^{l-1})}{1-\rho_{3s}}\\
&+\sqrt{2(1+\theta_{3s})})\|e\|_{2}\\
&\leq\rho_{3s}^{l}\|\widehat{x}_{\{p+1,\ldots,s\}}\|_{2}\\
&+(\sqrt{2}\gamma_{3s}\rho_{3s}^{l-1}\frac{\sqrt{1+\delta_{3s}}}{1-\delta_{3s}}+\frac{\sqrt{2}\gamma_{3s}\tau_{3s}(1-\rho_{3s}^{l-1})}{1-\rho_{3s}}\\
&+\sqrt{2(1+\theta_{3s})})\|e\|_{2}.\\
&\leq\rho_{3s}^{l}\|\widehat{x}_{\{p+1,\ldots,s\}}\|_{2}+\omega_{3s}\|e\|_{2}.\\
\end{aligned}
\end{equation*}

\begin{theorem}\label{www}
Suppose the measurement $y=Ax$ with an $s$-sparse vector $x$. If the P-RIP and RIP constants of $A$ satisfies $\rho_{3s}^{2}<1$, then the $s$-sparse vector $x$ is recovered form the measurement vector $y$ via a number $n$ of iterations of NST+HT+FB satisfying
\begin{equation*}
n\leq\frac{\ln(2/\rho_{3s})}{\ln(1/\rho_{3s})}s.
\end{equation*}
\end{theorem}

$\bm{Proof.}$
Let $\pi$ be the permutation of $\{1,2,\ldots,N\}$ such that $|x_{\pi(i)}|=\widehat{x}_{i}$, $i\in\{1,2,\ldots,N\}$. The goal is to determine an integer $n$ so that \supp$(x)\subseteq T_{n}$. We make a partition $Q_{1}\cup Q_{2}\cup\ldots Q_{r}$, $r\leq s$ of \supp$(x)=\pi(\{1,\ldots,s\})$. The sets are defined as follows
\begin{equation}\label{han1}
Q_{i}=\pi(\{q_{i-1}+1,\ldots,q_{i}\}),
\end{equation}
where $q_{0}=0, q_{i}=$maximum index~$\geq q_{i-1}+1$ ~so that~$\widehat{x}_{q_{i}}>\frac{\widehat{x}_{q_{i-1}+1}}{\sqrt{2}}$.
By the definition, $\widehat{x}_{q_{i}+1}\leq \frac{\widehat{x}_{q_{i-1}+1}}{\sqrt{2}}$ for all $i\in\{1,\ldots,r-1\}.$

\noindent With the introduction of $Q_{0}=\emptyset$, $k_{0}=0$ and the definition
\begin{equation}\label{han2}
k_{i}=\lceil\frac{\ln(2(|Q_{i}|+|Q_{i+1}|/2\cdots+|Q_{r}|/2^{r-i}))}{\ln(1/\rho_{3s}^{2})}\rceil,
\end{equation}
we prove by induction on $i\in\{0,\ldots,r\}$ that
\begin{equation}\label{han3}
Q_{0}\cup Q_{1}\cup\ldots Q_{i}\subseteq T_{k_{0}+k_{1}+\ldots k_{i}}.
\end{equation}

\noindent For $i=0$, (\ref{han3}) holds trivially. Then if (\ref{han3}) holds for $i-1,i\in\{1£¬\ldots£¬r\}$, {\em Lemma \ref{qw}} guarantees that (\ref{han3}) holds for $i$, provided
\begin{equation}\label{han4}
(\widehat{x}_{q_{i}})^{2}>\rho_{3s}^{2k_{i}}(\|x_{Q_{i}}\|_{2}^{2}+\|x_{Q_{i+1}}\|_{2}^{2}+\cdots+\|x_{Q_{r}}\|_{2}^{2}).
\end{equation}

\noindent Due to $(\widehat{x}_{q_{i}+1})\leq\frac{(\widehat{x}_{q_{i-1}+1})}{\sqrt{2}},$
\begin{equation*}
\begin{aligned}
&\|x_{Q_{i}}\|_{2}^{2}+\|x_{Q_{i+1}}\|_{2}^{2}+\cdots+\|x_{Q_{r}}\|_{2}^{2}\\
&\leq(\widehat{x}_{q_{i-1}+1})^{2}|Q_{i}|+(\widehat{x}_{q_{i}+1})^{2}|Q_{i+1}|+\cdots+(\widehat{x}_{q_{r-1}+1})^{2}|Q_{r}|\\
&\leq(\widehat{x}_{q_{i-1}+1})^{2}(|Q_{i}|+\frac{1}{2}|Q_{i+1}|+\cdots+\frac{1}{2^{r-i}}|Q_{r}|)\\
&\leq2(\widehat{x}_{q_{i}})^{2}(|Q_{i}|+\frac{1}{2}|Q_{i+1}|+\cdots+\frac{1}{2^{r-i}}|Q_{r}|).\\
\end{aligned}
\end{equation*}
\noindent The condition (\ref{han4}) holds because of the definition of $k_{i}$ (\ref{han2}).
The inductive proof (\ref{han3}) is concluded. The number of iterations that derives the support of $x$ can be represented as
\begin{equation*}
\begin{aligned}
&n=\sum_{i=1}^{r}k_{i}\leq\sum_{i=1}^{r}(1+\frac{\ln(2(|Q_{i}|+|Q_{i+1}|/2\cdots+|Q_{r}|/2^{r-i}))}{\ln(1/\rho_{3s}^{2})})\\
&=r+\frac{r}{\ln(1/\rho_{3s}^{2})}\sum_{i=1}^{r}\frac{1}{r}\ln(2(|Q_{i}|+|Q_{i+1}|/2\cdots+|Q_{r}|/2^{r-i})).\\
\end{aligned}
\end{equation*}
\noindent With the concavity of the function $\ln(x)$,
\begin{equation*}
\begin{aligned}
\ln(1/\rho_{3s}^{2})\frac{n-r}{r}&\leq\sum_{i=1}^{r}\frac{1}{r}\ln(2(|Q_{i}|+|Q_{i+1}|/2\cdots+|Q_{r}|/2^{r-i}))\\
&\leq\ln(\sum_{i=1}^{r}\frac{2}{r}(|Q_{i}|+|Q_{i+1}|/2\cdots+|Q_{r}|/2^{r-i})\\
&=\ln(\frac{2}{r}(|Q_{1}|+(1+\frac{1}{2})|Q_{2}|+\cdots\\
&+(1+\frac{1}{2}+\cdots+\frac{1}{2^{r-1}})|Q_{r}|))\\
&\leq\ln(\frac{4}{r}(|Q_{1}|+|Q_{2}|+\cdots+|Q_{r}|))\\
&=\ln(\frac{4s}{r}).
\end{aligned}
\end{equation*}
After simplification, we have
\begin{equation*}
n\leq r+\frac{r\ln(\frac{4s}{r})}{\ln(1/\rho_{3s}^{2})}=\frac{r}{s}(s+\frac{s\ln(\frac{4s}{r})}{\ln(1/\rho_{3s}^{2})}).
\end{equation*}
\noindent Since the function $f(x)=\frac{1}{x}(s+\frac{s\ln(4x)}{\ln(1/\rho_{3s})})$ is monotone decreasing function when $x\geq1$, we derive that
\begin{equation}\label{han5}
n\leq (s+\frac{s\ln(4)}{\ln(1/\rho_{3s}^{2})})=\frac{\ln(2/\rho_{3s})}{\ln(1/\rho_{3s})}s.
\end{equation}

\begin{remark}
For instance, if the P-RIP and RIP constants of $A$ satisfies $\rho_{3s}<\frac{1}{2}$, then the $s$-sparse vector $x$ is recovered form the measurement vector $y$ via a number $n$ of iterations of NST+HT+FB satisfying $n\leq2s$. Furthermore, if $A$ is Parseval frame, $\delta_{3s}\leq\frac{\sqrt{3}}{6}$ yields $n\leq2s$.
\end{remark}

The above arguments probe the idealized situation, the following lemmas examine the realistic situation that the measurement is perturbed by additive noise. Compared with the smallest nonzero absolute entry of the sparse signal, the noise is not too large. Under the assumption, the sparse signal can be recovered in a number of iterations, independently of the types of noise and sparse signals.

\begin{theorem}
Suppose the measurement $y=Ax+e$ with an $s$-sparse vector $x$ and $\|e\|_{2}\leq\frac{(1/\sqrt{2}-1/4)\widehat{x}_{s}}{\omega_{3s}}$. If the P-RIP and RIP constants of $A$ satisfies $\rho_{3s}^{2}<\frac{\sqrt{2}}{8}$, then the $s$-sparse vector $x$ is recovered form the measurement vector $y$ via a number $n$ of iterations of NST+HT+FB satisfying
\begin{equation*}
\|x-\mu^{n}\|_{2}\leq\frac{\|e\|_{2}}{\sqrt{1-\delta_{3s}}},
\end{equation*}
where $n\leq 3s.$
\end{theorem}
$\bm{Proof.}$
 Let $\pi$ be the permutation of $\{1,2,\ldots,N\}$ such that $|x_{\pi(i)}|=\widehat{x}_{i}$, $i\in\{1,2,\ldots,N\}$. The goal of the proof is to determine an integer $n$ so that \supp$(x)\subseteq T_{n}$. We partition \supp$(x)=\pi(\{1,\ldots,s\})$ as $Q_{1}\cup Q_{2}\cup\ldots Q_{r}$, $r\leq s$, where $Q_{i}$ are defined in
(\ref{han1}). With $Q_{0}=\emptyset$, $k_{0}=0$, and the definition
\begin{equation}\label{han7}
k_{i}=\lceil\frac{\ln(16(|Q_{i}|+|Q_{i+1}|/2\cdots+|Q_{r}|/2^{r-i}))}{\ln(1/\rho_{3s}^{2})}\rceil,
\end{equation}
\noindent we prove by induction on $i\in\{0,\ldots,r\}$ that
\begin{equation}\label{han8}
Q_{0}\cup Q_{1}\cup\ldots Q_{i}\subseteq T_{k_{0}+k_{1}+\ldots k_{i}}.
\end{equation}

\noindent For $i=0$, (\ref{han8}) holds trivially. Then if (\ref{han8}) holds for $i-1, i\in\{1£¬\ldots£¬r\}$, {\em Lemma \ref{qw}} guarantees that (\ref{han8}) holds for $i$, provided
\begin{equation}\label{han9}
\widehat{x}_{q_{i}}>\rho_{3s}^{k_{i}}\sqrt{\|x_{Q_{i}}\|_{2}^{2}+\|x_{Q_{i+1}}\|_{2}^{2}+\cdots+\|x_{Q_{r}}\|_{2}^{2}}+\omega_{3s}\|e\|_{2}.
\end{equation}

\noindent Since $\widehat{x}_{q_{i}+1}\leq\frac{(\widehat{x}_{q_{i-1}+1})}{\sqrt{2}}$ and $\widehat{x}_{q_{i}}>\frac{\widehat{x}_{q_{i-1}+1}}{\sqrt{2}}$
\begin{equation*}
\begin{aligned}
&\sqrt{\|x_{Q_{i}}\|_{2}^{2}+\|x_{Q_{i+1}}\|_{2}^{2}+\cdots+\|x_{Q_{r}}\|_{2}^{2}}\\
&\leq\sqrt{(\widehat{x}_{q_{i-1}+1})^{2}|Q_{i}|+(\widehat{x}_{q_{i}+1})^{2}|Q_{i+1}|+\cdots+(\widehat{x}_{q_{r-1}+1})^{2}|Q_{r}|}\\
&\leq \widehat{x}_{q_{i-1}+1}\sqrt{|Q_{i}|+\frac{1}{2}|Q_{i+1}|+\cdots+\frac{1}{2^{r-i}}|Q_{r}|}\\
&\leq \sqrt{2}\widehat{x}_{q_{i}}\sqrt{|Q_{i}|+\frac{1}{2}|Q_{i+1}|+\cdots+\frac{1}{2^{r-i}}|Q_{r}|}.\\
\end{aligned}
\end{equation*}
\noindent The condition (\ref{han9}) holds because of the definition of $k_{i}$ (\ref{han7}), $\|e\|_{2}\leq\frac{(1/\sqrt{2}-1/4)\widehat{x}_{s}}{\omega_{3s}}$ and $\widehat{x}_{s}\leq \widehat{x}_{q_{i-1}+1}$.
The inductive proof (\ref{han8}) is concluded. The number of iterations that derives the support of $x$ can be represented as

\begin{equation*}\label{han10}
n\leq (s+\frac{s\ln(32)}{\ln(1/\rho_{3s}^{2})})\leq3s.
\end{equation*}
When $S\subset T_{n}$, the feedback from $x^{n}$ to $\mu^{n}$ is equal to zero. Due to $x^{n}$ is a feasible solution, $\mu^{n}$ is also a feasible solution, i.e., $y=A\mu^{n}$.
\begin{equation*}\label{han11}
\begin{aligned}
\|x-\mu^{n}\|_{2}&\leq\frac{1}{\sqrt{1-\delta_{n}}}\|A(x-\mu^{n})\|_{2}\leq\frac{1}{\sqrt{1-\delta_{n}}}\|y-e-A\mu^{n}\|_{2}\\
&\leq\frac{\|e\|_{2}}{\sqrt{1-\delta_{3s}}}.
\end{aligned}
\end{equation*}

\begin{remark}
 Furthermore, if $A$ is the Parseval frame, the condition is relaxed to the RIP condition, i.e., $\delta_{3s}\leq\sqrt{\frac{1}{1+8\sqrt{2}}}$. The $s$-sparse vector $x$ is recovered form the measurement vector $y$ via a number $n$ of iterations of NST+HT+FB satisfying
\begin{equation*}
\|x-\mu^{n}\|_{2}\leq\frac{1}{\sqrt{1-1/\sqrt{1+8\sqrt{2}}}}\|e\|_{2},
\end{equation*}
where $n\leq3s.$
\end{remark}
\subsection{Specific signals recovery}
As shown in {\em Theorem \ref{www}}, NST+HT+FB can recover an $s$-sparse signals via the number of iterations $cs$, where $c>1$. In fact, for some specific $s$-sparse signals, the number of iterations can be significantly lowered by exploiting the
structures of signals. The following arguments correspond to the sparse Gaussian signal, which can be recovered via a number of iterations at most proportional to $\ln(s)$. For verifying the arguments, we first give the following lemma.
\begin{lemma}\label{propos11}
Suppose that the matrix $A$ satisfies
\begin{equation*}
\begin{aligned}
\delta_{2s}^{2}+2\gamma_{3s}^{2}<1,
\end{aligned}
\end{equation*}
then any $s$-sparse signal $x$ can be recovered by NST+HT+FB with $y=Ax$ in at most
\begin{equation}\label{start}
\begin{aligned}
\lceil2\ln(\frac{\|x\|_{2}}{\widehat{x}_{s}})/\ln(\frac{1-\delta_{2s}^{2}}{2\gamma_{3s}^{2}})\rceil
\end{aligned}
\end{equation}
iterations.
\end{lemma}
$\bm{Proof.}$
According to the stopping condition, the proof needs to determine an integer $k$ such that $T_{k}=S$. It is enough to show that, for all $i\in S$ and all $j\in S^{c}$,
\begin{equation}\label{dddd}
\begin{aligned}
&|[\mu^{k-1}+A^{\ast}(AA^{\ast})^{-1}(y-A\mu^{k-1})]_{i}|\\
&>|[\mu^{k-1}+A^{\ast}(AA^{\ast})^{-1}(y-A\mu^{k-1})]_{j}|.\\
\end{aligned}
\end{equation}

We notice that
\begin{equation*}
\begin{aligned}
&|[\mu^{k-1}+A^{\ast}(AA^{\ast})^{-1}(y-A\mu^{k-1})]_{i}|\\
&=|x_{i}+[(I-A^{\ast}(AA^{\ast})^{-1}A)(\mu^{k-1}-x)]_{i}|\\
&\geq\widehat{x}_{s}-|[(I-A^{\ast}(AA^{\ast})^{-1}A)(\mu^{k-1}-x)]_{i}|\\
\end{aligned}
\end{equation*}
and
\begin{equation*}
\begin{aligned}
&|[\mu^{k-1}+A^{\ast}(AA^{\ast})^{-1}(y-A\mu^{k-1})]_{j}|\\
&=|[(I-A^{\ast}(AA^{\ast})^{-1}A)(\mu^{k-1}-x)]_{j}|.\\
\end{aligned}
\end{equation*}
Using {\em Remark \ref{hhhh1}}, one can derive
\begin{equation*}
\begin{aligned}
&|[(I-A^{\ast}(AA^{\ast})^{-1}A)(\mu^{k-1}-x)]_{i}|\\
&+|[(I-A^{\ast}(AA^{\ast})^{-1}A)(\mu^{k-1}-x)]_{j}|\\
&\leq\sqrt{2}\|[(I-A^{\ast}(AA^{\ast})^{-1}A)(\mu^{k-1}-x)]_{\{i,j\}}\|_{2}\\
&\leq\sqrt{2}\gamma_{3s}\|\mu^{k-1}-x\|_{2}\\
&\leq(\frac{2\gamma_{3s}^{2}}{(1-\delta_{2s}^{2})})^{\frac{k}{2}}\|\mu^{0}-x\|_{2}.
\end{aligned}
\end{equation*}
Therefore, (\ref{dddd}) is fulfilled as soon as
\begin{equation*}
\begin{aligned}
\widehat{x}_{s}\geq(\frac{2\gamma_{3s}^{2}}{(1-\delta_{2s}^{2})})^{\frac{k}{2}}\|\mu^{0}-x\|_{2},
\end{aligned}
\end{equation*}
i.e.,
\begin{equation*}
\begin{aligned}
k\geq2\ln(\frac{\|x\|_{2}}{\widehat{x}_{s}})/\ln(\frac{1-\delta_{2s}^{2}}{2\gamma_{3s}^{2}}),
\end{aligned}
\end{equation*}
where $\mu^{0}=0$.
The smallest such integer $k$ is given by (\ref{start}).

\begin{proposition}
Suppose that the matrix $A$ satisfies
$\delta_{2s}^{2}+2\gamma_{3s}^{2}<1$
and an $s$-sparse signal $x$ is taken as independent standard Gaussian random variables $x_{i}$, where $i\in$supp$(x)$. Then, with probability larger $1-\eta$,
the $s$-sparse signal $x$ is recovered from $y=Ax$ via a number of NST+HT+FB iterations at most proportional to $\ln(\frac{s}{\eta})$.
\end{proposition}
$\bm{Proof.}$
Since $\widehat{x}_{i}$ follows the standard Gaussian distribution, for $i\in$supp$(x)$, then
\begin{equation}\label{ddddd}
\begin{aligned}
\mathbb{P}(\widehat{x}_{s}<t)&=\mathbb{P}(|x_{i}|<t, i\in supp(x))\leq s\mathbb{P}(|x|<t)\\
&\leq s\int_{-t}^{t}\frac{\exp(-\nu^{2}/2)}{\sqrt{2\pi}}d\nu\\
&\leq s\int_{-t}^{t}\frac{1}{\sqrt{2\pi}}d\nu=s\sqrt{\frac{2}{\pi}}t,
\end{aligned}
\end{equation}
\begin{equation}\label{sssss}
\begin{aligned}
&\mathbb{P}(\|x\|_{2}^{2}>t)\\
&=\mathbb{P}(\sum_{i\in S}x_{i}^{2}>t)\leq\frac{\mathbb{E}(\exp(\frac{\sum_{i\in S}x_{i}^{2}}{4}))}{\exp(\frac{t}{4})}=\frac{\prod_{i\in S}\mathbb{E}(\exp(\frac{x_{i}^{2}}{4}))}{\exp(\frac{t}{4})}\\
&=\frac{\prod_{i\in S}\int\exp(\frac{x_{i}^{2}}{4})\frac{1}{\sqrt{2\pi}}\exp(-\frac{x_{i}^{2}}{2})dx_{i}}{\exp(\frac{t}{4})}=\frac{2^{\frac{s}{2}}}{\exp(\frac{t}{4})}\\
&\leq\exp(\frac{2s-t}{4}),\\
\end{aligned}
\end{equation}
where the first inequality in (\ref{sssss}) is based on the $Markov's$ inequality.
Let $t=\sqrt{\frac{\pi}{8}}\frac{\eta}{s}$ in (\ref{ddddd}) and $t=2s-4\ln(\frac{\eta}{2})$ in (\ref{sssss}), where $\eta<1$, we have
\begin{equation*}
\begin{aligned}
&\mathbb{P}(\widehat{x}_{s}<\sqrt{\frac{\pi}{8}}\frac{\eta}{s})\leq\frac{\eta}{2},\\
&\mathbb{P}(\|x\|_{2}^{2}>2s-4\ln(\frac{\eta}{2}))\leq\frac{\eta}{2}.
\end{aligned}
\end{equation*}
Furthermore, let $\beta>2\eta-\frac{4\eta}{s}\ln(\frac{\eta}{2})$, then $\mathbb{P}(\|x\|_{2}^{2}>\beta\frac{s}{\eta})\leq\frac{\eta}{2}$.
Therefore, with probability larger than $1-\eta$, $\widehat{x}_{s}\geq\sqrt{\frac{\pi}{8}}\frac{\eta}{s}$ and $\|x\|_{2}^{2}\geq\beta\frac{s}{\eta}$. With {\em Lemma \ref{propos11}}, the number $k$ of iterations for recovering $x$ satisfies
\begin{equation*}
 k\leq\lceil2\ln(\sqrt{\beta}\sqrt{\frac{8}{\pi}}(\frac{s}{\eta})^{\frac{3}{2}})/\ln(\frac{1-\delta_{2s}^{2}}{2\gamma_{3s}^{2}})\rceil.
\end{equation*}

\section{The number of iterations of AdptNST+HT+FB}
In this section, the theoretical analysis of AdptNST+HT+FB is presented. The number of iterations for recovering an $s$-sparse signal is first established including the idealized situation and the realistic situation where the uniform sparse recovery can apply to all $s$-sparse signals. Then it is followed by the nonuniform setting.
\subsection{Uniform sparse recovery}
An analog of {\em Lemma \ref{qw}} can be obtained for AdptNST+HT+FB. The following Lemma shows that if the $p$ largest absolute entries are contained in the support at iteration $k$, then $l$ further iterations of AdptNST+HT+FB are sufficient to capture the indices of the $q$ following largest entries.
\begin{lemma}\label{tiu}
 Suppose $y=Ax+e$ with an $s$-sparse vector $x$ and let $\{T_{k}\}$ be index sets of the sequence $\{u^{k}\}$ in AdptNST+HT+FB. For integers $k\geq s$, $p\geq0$, if $T_{k}$ contains the indices of $p$ largest absolute entries of $x$. Then $\{T_{k+l}\}$ contains the indices of $p+q$ largest absolute entries of $x$ via $l$ iterations for integers $l,q\geq1$, provided
\begin{equation*}
\widehat{x}_{p+q}>\rho_{s+2k+2l}^{l}\|\widehat{x}_{\{p+1,\ldots,s\}}\|_{2}+\omega_{s+2k+2l}\|e\|_{2}.
\end{equation*}
\end{lemma}
$\bm{Proof.}$
For AdptNST+HT+FB, the hypothesis is $\pi(\{1,\ldots,p\})\subseteq T_{k}$ and the goal is to prove that $\pi(\{1,\ldots,p+q\})\subseteq T_{k+l}$. That is to say the $|(\mu^{l+k-1}+A^{\ast}(AA^{\ast})^{-1}(y-A\mu^{l+k-1}))_{\pi(j)}|$ for $j\in\{1,\ldots p+q\}$ are among the $l+k$ largest entries of
$|(\mu^{l+k-1}+A^{\ast}(AA^{\ast})^{-1}(y-A\mu^{l+k-1}))_{i}|$ for $i\in\{1,\ldots,N\}$. Since $l+k\geq s$, it is enough to prove that they are among the $s$ largest values of $|(\mu^{l+k-1}+A^{\ast}(AA^{\ast})^{-1}(y-A\mu^{l+k-1}))_{i}|$. The rest of the proof refers to the proof of {\em Lemma \ref{qw}}.

To make the algorithm more applicable, one may increase the sparsity by the iteration gradually. As discussed in {\em Lemma \ref{tiu}}, the cardinality of the intersection of $S$ and $T_{k}$ increases as the increasing of $k$. We can prove that when $S\subseteq T_{k}$ ($|T_{k}|=k$ and $k\geq s$), the AdptNST+HT+FB also converges.

\noindent$\bm{Proof.}$
Since  $S\subseteq T_{k}$,
\begin{equation*}
\begin{aligned}
 \eta^{k}&=\arg\min \limits_{\eta}\|A_{T_{k}}\eta-A_{T_{k}^{c}}x_{T_{k}^{c}}^{k}\|_{2}\\
 &=\arg\min \limits_{\eta}\|(A_{S}\eta_{S}+A_{T_{k}\setminus S}\eta_{T_{k}\setminus S})-(y-A_{T_{k}}x_{T_{k}}^{k})\|_{2}\\
 &=\arg\min \limits_{\eta}\|(A_{S}\eta_{S}+A_{S}x_{S}^{k})+(A_{T_{k}\setminus S}\eta_{T_{k}\setminus S}+A_{T_{k}\setminus S}x_{T_{k}\setminus S}^{k})\\
 &-y\|_{2}.
\end{aligned}
\end{equation*}
$\eta_{S}^{k}=x_{S}-x_{S}^{k}, \eta_{T_{k}\setminus S}^{k}=-x_{T_{k}\setminus S}^{k}$ is an solution. Then we can obtain $A_{T_{k}}\eta^{k}-A_{T_{k}^{c}}x_{T_{k}^{c}}^{k}=0$ so that $y=Au^{k}$. Due to the NST step $x^{k+1}=u^{k}+A^{\ast}(AA^{\ast})^{-1}(y-Au^{k})$ and the processing of feedback, $\mu^{k+1}=\mu^{k}$.  Therefore, $S\subseteq T_{k}$  is also the stopping criteria for AdptNST+HT+FB. The remaining topic is to determine the smallest integer $k$ such that $S\subseteq T_{k}$.

\begin{theorem}\label{Te3}
 Suppose the measurement $y=Ax$ with an $s$-sparse vector $x$. If the P-RIP and RIP constants of $A$ satisfies $\rho_{5s+\lceil\frac{\ln4}{\ln\alpha}\rceil s}^{2}<\frac{1}{\alpha^{2}}$, where $\alpha\in(1,+\infty),$ then the $s$-sparse vector $x$ is recovered form the measurement vector $y$ via a number $n$ of iterations of AdptNST+HT+FB satisfying
\begin{equation*}
n\leq(\frac{\ln2}{\ln\alpha}+2)s.
\end{equation*}
\end{theorem}

$\bm{Proof.}$
 Let $\pi$ be the permutation of $\{1,2,\ldots,N\}$ such that $|x_{\pi_{i}}|=\widehat{x}_{i}$ for all $i\in\{1,2,\ldots,N\}$. As discussed, $S\subseteq T_{k}(k\geq s)$ is the stopping condition for AdptNST+HT+FB. The goal is to determine an integer $n$ so that $S\subseteq T_{n}$. We still partition \supp$(x)=\pi(\{1,\ldots,s\})$ as $Q_{1}\cup Q_{2}\cup\ldots Q_{r}$, $r\leq s$, where $Q_{i}$ are defined in (\ref{han1}). Since $|T_{k}|\leq s $ in the first $s$ iteration, the first $s$ iteration are ignored. With $Q_{0}=\emptyset$ and $k_{0}=0$,
\begin{equation}\label{han12}
k_{i}=\lceil\frac{\ln(2(|Q_{i}|+|Q_{i+1}|/2\cdots+|Q_{r}|/2^{r-i}))}{\ln(1/\rho_{5s+\lceil\frac{\ln4}{\ln\alpha}\rceil s}^{2})}\rceil,
\end{equation}
we prove by induction on $i\in\{0,\ldots,r\}$ that
\begin{equation}\label{han13}
Q_{0}\cup Q_{1}\cup\ldots Q_{i}\subseteq T_{s+k_{0}+k_{1}+\ldots k_{i}}.
\end{equation}
For $i=0$, (\ref{han13}) holds trivially. Then if (\ref{han13}) holds for $i-1,~i\in\{1£¬\ldots£¬r\}$, {\em Lemma \ref{tiu}} guarantees that (\ref{han13}) holds for $i$, provided
\begin{equation} \label{han14}
(\widehat{x}_{q_{i}})^{2}>\rho_{s+2k_{i}+2(s+k_{1}+\cdots+k_{i-1})}^{2k_{i}}(\|x_{Q_{i}}\|_{2}^{2}+\|x_{Q_{i+1}}\|_{2}^{2}+\cdots+\|x_{Q_{r}}\|_{2}^{2}).
\end{equation}
As the same proof of (\ref{han4}), we have
\begin{equation}\label{han15}
(\widehat{x}_{q_{i}})^{2}>\rho_{5s+\lceil\frac{\ln4}{\ln\alpha}\rceil s}^{2k_{i}}(\|x_{Q_{i}}\|_{2}^{2}+\|x_{Q_{i+1}}\|_{2}^{2})+\cdots+\|x_{Q_{r}}\|_{2}^{2}).
\end{equation}
\noindent Through the definition of $k_{i}$ (\ref{han12}) and the inference of (\ref{han5}), we have
\begin{equation*}
\sum_{i=1}^{r}k_{i}\leq\frac{\ln(4/\rho_{5s+\lceil\frac{\ln4}{\ln\alpha}\rceil s}^{2})}{\ln(1/\rho_{5s+\lceil\frac{\ln4}{\ln\alpha}\rceil s}^{2})}s\leq(\frac{\ln2}{\ln\alpha}+1)s.
\end{equation*}
 \noindent Therefore, $\rho_{s+2k_{i}+2(s+k_{1}+\cdots+k_{i-1})}\leq\rho_{5s+\lceil\frac{\ln4}{\ln\alpha}\rceil s}$. Since (\ref{han15}) holds, the condition (\ref{han14}) is obviously fulfilled.
 The number of iterations that derives the support of $x$ can be represented as
\begin{equation}\label{han16}
n\leq\frac{\ln(4/\rho_{5s+\lceil\frac{\ln4}{\ln\alpha}\rceil s}^{2})}{\ln(1/\rho_{5s+\lceil\frac{\ln4}{\ln\alpha}\rceil s}^{2})}s+s\leq(\frac{\ln2}{\ln\alpha}+2)s.
\end{equation}

\begin{remark}
{\em Theorem \ref{Te3}} shows $\{u^{k}\}$ in AdptNST+HT+FB converges to $x$ in finitely many steps. We can observe that the condition of RIP and P-RIP becomes weaker and the upper bound of $n$ increases as $\alpha$ decreases. For instance, $\alpha=2$ yields the RIP and P-RIP condition $\rho_{7s}^{2}<\frac{1}{4}$ and $n\leq3s$, while $\alpha=\sqrt{2}$ yields the weaker RIP and P-RIP condition $\rho_{9s}^{2}<\frac{1}{2}$ and $n\leq4s$. Furthermore, if $A$ is the Parseval frame, $\alpha=2$ yields the RIP condition $\delta_{7s}<\frac{\sqrt{1}}{3},$ and $n\leq3s$, while $\alpha=\sqrt{2}$ yields the weaker RIP condition $\delta_{9s}<\frac{\sqrt{5}}{5},$ and $n\leq4s$.
\end{remark}

With the idealized situation in {\em Theorem \ref{Te3}}, we then extend the noiseless sampled data to the noisy case. The result shows that an $s$-sparse signal can be recovered under the RIP and P-RIP condition in the realistic situation and the error bound depend on the noise.

\begin{theorem}\label{Te5}
Suppose the measurement $y=Ax+e$ with an $s$-sparse vector $x$ and $\|e\|_{2}\leq\frac{(1/\sqrt{2}-1/4)\widehat{x}_{s}}{\omega_{5s+\lceil\frac{\ln32}{\ln\alpha}\rceil s}}$.  If the P-RIP and RIP constant of $A$ satisfies $\rho_{5s+\lceil\frac{\ln32}{\ln\alpha}\rceil s}^{2}<\frac{1}{\alpha^{2}}$, where $\alpha\in(1,+\infty),$ then the $s$-sparse vector $x$ is recovered from the measurement vector $y$ via a number $n$ of iterations of AdptNST+HT+FB satisfying
\begin{equation*}
\|x-\mu^{n}\|_{2}\leq\frac{\|e\|_{2}}{\sqrt{1-\delta_{2s+\lceil\frac{\ln32}{\ln\alpha^{2}})\rceil s}}},
\end{equation*}
where $n\leq (2+\frac{\ln32}{\ln\alpha^{2}})s.$
\end{theorem}

$\bm{Proof.}$
Let $\pi$ be the permutation of $\{1,2,\ldots,N\}$, i.e, $|x_{\pi_{(i)}}|=\widehat{x}_{i}$ for $i\in\{1,2,\ldots,N\}$. The goal of the proof is to determine an integer $n$ so that \supp$(x)\subseteq T_{n}$. We partition \supp$(x)=\pi(\{1,\ldots,s\})$ as $Q_{1}\cup Q_{2}\cup\ldots Q_{r}$, $r\leq s$, where $Q_{i}$ are defined in
(\ref{han1}). Since $|T_{k}|\leq s $ in the first $s$ iteration, the first $s$ iteration are ignored. With $Q_{0}=\emptyset$ and $k_{0}=0$,
\begin{equation}\label{han20}
k_{i}=\lceil\frac{\ln(16(|Q_{i}|+|Q_{i+1}|/2\cdots+|Q_{r}|/2^{r-i}))}{\ln(1/\rho_{5s+\lceil\frac{\ln32}{\ln\alpha}\rceil s}^{2})}\rceil,
\end{equation}
we prove by induction on $i\in\{0,\ldots,r\}$ that
\begin{equation}\label{han21}
Q_{0}\cup Q_{1}\cup\ldots Q_{i}\subseteq T_{s+k_{0}+k_{1}+\ldots k_{i}}.
\end{equation}
For $i=0$, (\ref{han21}) holds trivially. Then if (\ref{han21}) holds for $i-1$, $i\in\{1£¬\ldots£¬r\}$,  {\em Lemma \ref{tiu}} guarantees that (\ref{han21}) holds for $i$, provided
\begin{equation}\label{han22}
\begin{aligned}
&(\widehat{x}_{q_{i}})>\\
&\rho_{s+2k_{i}+2(s+k_{1}+\cdots+k_{i-1})}^{k_{i}}\sqrt{\|x_{Q_{i}}\|_{2}^{2}+\|x_{Q_{i+1}}\|_{2}^{2}+\cdots+\|x_{Q_{r}}\|_{2}^{2}}\\
&+\omega_{s+2k_{i}+2(s+k_{1}+\cdots+k_{i-1})}\|e\|_{2}.\\
\end{aligned}
\end{equation}
\noindent As the same proof of (\ref{han9}), we have
\begin{equation}\label{han23}
\begin{aligned}
&\widehat{x}_{q_{i}}>\\
&\rho_{5s+\lceil\frac{\ln32}{\ln\alpha}\rceil s}^{k_{i}}\sqrt{\|x_{Q_{i}}\|_{2}^{2}+\|x_{Q_{i+1}}\|_{2}^{2}+\cdots+\|x_{Q_{r}}\|_{2}^{2}}\\&
+\omega_{5s+\lceil\frac{\ln32}{\ln\alpha}\rceil s}\|e\|_{2}
\end{aligned}
\end{equation}

\noindent Through the definition of $k_{i}$ (\ref{han20}) and the inference of (\ref{han5}), we have
\begin{equation*}
\sum_{i=1}^{r}k_{i}\leq\frac{\ln(32/\rho_{5s+\lceil\frac{\ln32}{\ln\alpha}\rceil s}^{2})}{\ln(1/\rho_{5s+\lceil\frac{\ln32}{\ln\alpha}\rceil s}^{2})}s\leq(1+\frac{\ln32}{\ln\alpha^{2}})s.
\end{equation*}
\noindent Therefore, $\rho_{s+2k_{i}+2(s+k_{1}+\cdots+k_{i-1})}\leq\rho_{5s+\lceil\frac{\ln32}{\ln\alpha}\rceil s}$ and $\omega_{s+2k_{i}+2(s+k_{1}+\cdots+k_{i-1})}\leq\omega_{5s+\lceil\frac{\ln32}{\ln\alpha}\rceil s}$. Since (\ref{han22}) holds, the condition
(\ref{han21}) is obviously fulfilled.

 The number of iterations that derives the support of $x$ can be represented as
\begin{equation*}
n\leq\frac{\ln(32/\rho_{5s+\lceil\frac{\ln32}{\ln\alpha}\rceil s}^{2})}{\ln(1/\rho_{5s+\lceil\frac{\ln32}{\ln\alpha}\rceil s}^{2})}s\leq(2+\frac{\ln32}{\ln\alpha^{2}})s.
\end{equation*}
It then follows that
\begin{equation*}
\begin{aligned}
\|x-\mu^{n}\|_{2}&\leq\frac{1}{\sqrt{1-\delta_{n}}}\|A(x-\mu^{n})\|_{2}\\
&\leq\frac{1}{\sqrt{1-\delta_{n}}}\|y-e-A\mu^{n}\|_{2}\\
&\leq\frac{\|e\|_{2}}{\sqrt{1-\delta_{2s+\lceil\frac{\ln32}{\ln\alpha^{2}}\rceil s}}}.\\
\end{aligned}
\end{equation*}

\begin{remark}
{\em Theorem \ref{Te5}} shows that the condition of RIP and P-RIP becomes weaker and the upper bound of $n$ and error increases with $\alpha$ decreasing. For instance, $\alpha=4$ yields the RIP and P-RIP condition $\rho_{8s}^{2}<\frac{1}{16}$, $n\leq4s$, $\|x-\mu^{n}\|_{2}\leq\frac{\|e\|_{2}}{\sqrt{1-\delta_{4s}}},$ while
$\alpha=2$ yields the  weaker RIP and P-RIP condition $\rho_{10s}^{2}<\frac{1}{4}$, $n\leq5s$, and $\|x-\mu^{n}\|_{2}\leq\frac{\|e\|_{2}}{\sqrt{1-\delta_{5s}}}.$ Furthermore, if $A$ is A is Parseval frame, $\alpha=4$ yields the RIP condition $\delta_{8s}<\frac{\sqrt{33}}{33},$ $n\leq4s$, $\|x-\mu^{n}\|_{2}\leq\frac{\sqrt{33}\|e\|_{2}}{\sqrt{33-\sqrt{33}}},$ while
$\alpha=2$ yields the weaker RIP condition $\delta_{10s}<\frac{1}{3},$ $n\leq5s$, and $\|x-\mu^{n}\|_{2}\leq\frac{\sqrt{6}\|e\|_{2}}{2}.$
\end{remark}
\subsection{Specific signals recovery}
{\em Theorem \ref{Te3}} and  {\em Theorem \ref{Te5}} show that the number of iterations of AdptNST+HT+FB is greater than $2s$ for recovering an $s$-sparse signal. In this section, we consider a specific case. The conclusion demonstrates that the recovery of an $s$-sparse signals via exactly $s$ iterations can be obtained in an specific setting. The proof of the conclusion follows the similar trajectory with \cite{R27,R28}.
The following results play an important role in the proof of propositions.
\begin{lemma}\cite{R29}\label{Pro1}
Suppose that an $M\times N$ random matrix $A$ is drawn according to a probability distribution for which the concentration inequality holds, i.e., for $t\in(0,1)$, a constant $C_{1}\in\mathbb{R}$,
\begin{equation*}
\mathbb{P}(|\|Ax\|_{2}^{2}-\|x\|_{2}^{2}|\geq t\|x\|_{2}^{2})\leq2\exp(-C_{1}t^{2}M),~~for~all~x\in\mathbb{R}^{N}.
\end{equation*}
Then, for a fixed set $S\subset\{1,2,\cdots,N\}$ with cardinality $s$,
\begin{equation*}
\mathbb{P}(\|A_{S}^{\ast}A_{S}-I\|_{2\rightarrow2}>\delta)\leq2\exp(-C_{2}\delta^{2}M),
\end{equation*}
if $M\geq C_{3}s/\delta^{2}$, with $C_{2}$ and $C_{3}$ depending only the entries distributions.
\end{lemma}

\begin{remark}\label{rema1}
Suppose that an $M\times N$ random matrix $A$ is drawn according to a probability distribution and the matrix $B=((AA^{\ast})^{-\frac{1}{2}}A)$ satisfies the concentration inequality holds, i.e., for $t'\in(0,1)$, a constant $C_{1}'\in\mathbb{R}$,
\begin{equation*}
\mathbb{P}(|\|Bx\|_{2}^{2}-\|x\|_{2}^{2}|\geq t'\|x\|_{2}^{2})\leq2\exp(-C_{1}'(t')^{2}M),~~for~all~x\in\mathbb{R}^{N}.
\end{equation*}
Then, for a fixed set $S\subset\{1,2,\cdots,N\}$ with cardinality $s$,
\begin{equation*}
\mathbb{P}(\|B_{S}^{\ast}B_{S}-I\|_{2\rightarrow2}>\gamma)\leq2\exp(-C_{2}'\gamma^{2}M),
\end{equation*}
i.e.,
\begin{equation*}
\mathbb{P}(\|A_{S}^{\ast}(AA^{\ast})^{-1}A_{S}-I\|_{2\rightarrow2}>\gamma)\leq2\exp(-C_{2}'\gamma^{2}M),
\end{equation*}
if $M\geq C_{3}'s/\gamma^{2}$, with $C_{2}'$ and $C_{3}'$ depending only the entries distributions.
An analog of conclusions can be obtained for the matrix $E=((AA^{\ast})^{-1}A)$. For a fixed set $S\subset\{1,2,\cdots,N\}$ with cardinality $s$,
\begin{equation*}
\mathbb{P}(\|A_{S}^{\ast}(AA^{\ast})^{-2}A_{S}-I\|_{2\rightarrow2}>\theta)\leq2\exp(-C_{2}''\theta^{2}M),
\end{equation*}
if $M\geq C_{3}''s/\theta^{2}$, with $C_{2}''$ and $C_{3}''$ depending only the entries distributions.
\end{remark}

\begin{lemma}\cite{R29}\label{Pro2}
For an $M\times N$ random matrix $A$ with independent subgaussian entries and a vector $v$ and an index $l\in\{1,2,\cdots N\}$,
\begin{equation*}
\mathbb{P}(|\langle A_{l},v\rangle|\geq t\|v\|_{2})\leq4\exp(-C_{4}t^{2}M),
\end{equation*}
where the constant $C_{4}$ depends only the subgaussian distribution.
\end{lemma}

\begin{remark}\label{remark}
For a random matrix $B=((AA^{\ast})^{-\frac{1}{2}}A)$ defined in {\em Remark \ref{rema1}} and a vector $v$
\begin{equation*}
\mathbb{P}(|\langle B_{l},v\rangle|\geq t\|v\|_{2})\leq4\exp(-C_{4}'t^{2}M),
\end{equation*}
where the constant $C_{4}'$ depends only the distribution of $B$.
\end{remark}

\begin{proposition}\label{theorem1}
Suppose $y=Ax+e$ with an $s$-sparse signal $x$ such that $\widehat{x}_{1}\leq\sigma\widehat{x}_{s}$ and the noise such that $\|e\|_{2}\leq\epsilon\widehat{x}_{s}$, where $A\in\mathbb{R}^{M\times N}$ is a Gaussian random matrix and $\sigma\geq1$.
If the number of measurements satisfies
\begin{equation*}
M\geq Cs\ln(N),
\end{equation*}
then with probability larger than $1-6N^{-\alpha}$, the sequence of $\{\mu^{k}\}$ at iteration $s$ in AdptNST+HT+FB satisfies,
\begin{equation*}
T_{s}=\textit{supp}(x)~\textit{and}~\|x-\mu^{s}\|_{2}\leq p\|e\|_{2}.
\end{equation*}
The constants $\epsilon$ and $p$ depend only on $\sigma$, while the constant $C$ depends on $\epsilon$ and $\alpha$.
\end{proposition}
$\bm{Proof.}$
As defined, $\widehat{x}\in\mathbb{R}^{N}_{+}$ is the nonincreasing rearrangement of a vector $x=(x_1, x_2,\cdots,x_N)^{'}\in\mathbb{R}^{N}$, i.e., $\widehat{x}_{1}\geq \widehat{x}_{2}\geq\ldots \widehat{x}_{N}\geq0$
and there exists a permutation $\pi$ of $\{1,\ldots,N\}$ such that $\widehat{x}_{i}=|x_{\pi(i)}|$ for all $i\in\mathbb{R}^{N}$. We define two random variables $\psi_{k}$ and $\phi_{k}$ for all $k\in\{1,\ldots,s\}$ as
\begin{equation*}
\begin{aligned}
&\psi_{k}=\widehat{[(\mu^{k-1}+A^{\ast}(AA^{\ast})^{-1}(y-A\mu^{k-1}))_{S}]_{k}},\\
&\phi_{k}=\widehat{[(\mu^{k-1}+A^{\ast}(AA^{\ast})^{-1}(y-A\mu^{k-1}))_{S^{c}}]_{1}},\\
\end{aligned}
\end{equation*}
i.e., $\psi_{k}$ is the $k$th largest value of $|(\mu^{k-1}+A^{\ast}(AA^{\ast})^{-1}(y-A\mu^{k-1}))_{i}|$, $i\in S$ and $\phi_{k}$ is the largest value of $|(\mu^{k-1}+A^{\ast}(AA^{\ast})^{-1}(y-A\mu^{k-1}))_{j}|$, $j\in S^{c}$. It can be noted that $\psi_{k}>\phi_{k}$ for all $k\in\{1,\ldots,s\}$ indicates that $T_{k}\subseteq S$ for all $k\in\{1,\ldots,s\}$. The failure probability of this event is as follows
\begin{equation*}
\begin{aligned}
P&=\mathbb{P}(\exists~k\in\{1,\ldots,s\}:\phi_{k}\geq\psi_{k}\\
&~\text{and}~(\psi_{k-1}>\phi_{k-1},\ldots,\psi_{1}>\phi_{1}))\\
&\leq\sum_{j\in S^{c}}\mathbb{P}(\|A_{S\bigcup\{j\}}^{\ast}A_{S\bigcup\{j\}}-I\|_{2\rightarrow2}>\delta)\\
&+\sum_{j\in S^{c}}\mathbb{P}(\|A_{S\bigcup\{j\}}^{\ast}(AA^{\ast})^{-1}A_{S\bigcup\{j\}}-I\|_{2\rightarrow2}>\gamma)\\
&+\sum_{j\in S^{c}}\mathbb{P}(\|A_{S\bigcup\{j\}}^{\ast}(AA^{\ast})^{-2}A_{S\bigcup\{j\}}-I\|_{2\rightarrow2}>\theta)\\
&+\sum_{k=1}^{s}\mathbb{P}[\phi_{k}\geq\psi_{k},(\psi_{k-1}>\phi_{k-1},\ldots,\psi_{1}>\phi_{1}),\\
&(\|A_{S\bigcup\{j\}}^{\ast}A_{S\bigcup\{j\}}-I\|_{2\rightarrow2}\leq\delta ~\text{for}~\text{all}~j\in S^{c}),\\
&(\|A_{S\bigcup\{j\}}^{\ast}(AA^{\ast})^{-1}A_{S\bigcup\{j\}}-I\|_{2\rightarrow2}\leq\gamma ~\text{for}~\text{all}~j\in S^{c})),\\
&(\|A_{S\bigcup\{j\}}^{\ast}(AA^{\ast})^{-2}A_{S\bigcup\{j\}}-I\|_{2\rightarrow2}\leq\theta ~\text{for}~\text{all}~j\in S^{c}))].\\
\end{aligned}
\end{equation*}
According to {\em Lemma \ref{Pro1}} and {\em Remark \ref{rema1}}, the  first three terms of the inequality is bounded by
\begin{equation*}
\begin{aligned}
6(N-s)\exp(-\min\{C_{2},C_{2}',C_{2}''\}\min\{\delta^{2},\gamma^{2},\theta^{2}\}M)).
\end{aligned}
\end{equation*}
We now turn to the last term of the last inequality. For simplicity, we use $\overline{\mathbb{P}}(\phi_{k}\geq\psi_{k})$ to denote the probability of the event $\phi_{k}\geq\psi_{k}$ intersected with the events $(\psi_{k-1}>\phi_{k-1},\ldots,\psi_{1}>\phi_{1})$, $(\|A_{S\bigcup\{j\}}^{\ast}A_{S\bigcup\{j\}}-I\|_{2\rightarrow2}\leq\delta ~\text{for}~\text{all}~j\in S^{c})$, $(\|A_{S\bigcup\{j\}}^{\ast}(AA^{\ast})^{-1}A_{S\bigcup\{j\}}-I\|_{2\rightarrow2}\leq\gamma ~\text{for}~\text{all}~j\in S^{c}))$ and $(\|A_{S\bigcup\{j\}}^{\ast}(AA^{\ast})^{-2}A_{S\bigcup\{j\}}-I\|_{2\rightarrow2}\leq\theta ~\text{for}~\text{all}~j\in S^{c}))$. Let $T_{s-k+1}$ correspond to the support of the $s-k+1$ smallest values of $|(\mu^{k-1}+A^{\ast}(AA^{\ast})^{-1}(y-A\mu^{k}))_{S}|$. By the definition, we have
\begin{equation*}
\begin{aligned}
\psi_{k}&\geq\frac{1}{\sqrt{s-k+1}}\|(\mu^{k-1}+A^{\ast}(AA^{\ast})^{-1}(y-A\mu^{k-1}))_{T_{s-k+1}}\|_{2}\\
&\geq\frac{1}{\sqrt{s-k+1}}(\|x_{T_{s-k+1}}\|_{2}\\
&-\|((A^{\ast}(AA^{\ast})^{-1}A-I)(x-\mu^{k-1}))_{T_{s-k+1}}\|_{2}\\
&-\|((A^{\ast}(AA^{\ast})^{-1}e)_{T_{s-k+1}}\|_{2})\\.
\end{aligned}
\end{equation*}
Since $\psi_{k-1}>\phi_{k-1}$, we can obtain $T_{k-1}\subseteq S$. Therefore,
\begin{equation*}
\begin{aligned}
&\|((A^{\ast}(AA^{\ast})^{-1}A-I)(x-\mu^{k-1}))_{T_{s-k+1}}\|_{2}\\
&\leq\|(A_{S}^{\ast}(AA^{\ast})^{-1}A_{S}-I)(x-\mu^{k-1})\|_{2}\\
&\leq\gamma\|x-\mu^{k-1}\|_{2}.
\end{aligned}
\end{equation*}
By $\|((A^{\ast}(AA^{\ast})^{-1}e)_{T_{s-k+1}}\|_{2}\leq\sqrt{1+\theta}\|e\|_{2}$, where $\theta=\delta((AA^{\ast})^{-1}A)$.
It follows that
\begin{equation}\label{equa1}
\begin{aligned}
\psi_{k}\geq\frac{1}{\sqrt{s-k+1}}(\|x_{T_{s-k+1}}\|_{2}-\gamma\|x-\mu^{k-1}\|_{2}-\sqrt{1+\theta}\|e\|_{2}).
\end{aligned}
\end{equation}
Since $T_{k-1}\subseteq S$, one has $(\mu^{k-1})_{S^{c}}=0$. Then,
\begin{equation}\label{equa2}
\begin{aligned}
\phi_{k}&=\widehat{[(\mu^{k-1}+A^{\ast}(AA^{\ast})^{-1}(y-A\mu^{k-1}))_{S^{c}}]_{1}}\\
&=\widehat{[(A^{\ast}(AA^{\ast})^{-1}(y-A\mu^{k-1}))_{S^{c}}]_{1}}\\
&=\max_{j\in S^{c}}|(A^{\ast}(AA^{\ast})^{-1}(y-A\mu^{k-1}))_{j}|\\
&\leq\max_{j\in S^{c}}|(A^{\ast}(AA^{\ast})^{-1}A(x-\mu^{k-1}))_{j}|+\max_{j\in S^{c}}|(A^{\ast}(AA^{\ast})^{-1}e)_{j}|\\
&=\max_{j\in S^{c}}|[((AA^{\ast})^{-\frac{1}{2}}A)^{\ast}(AA^{\ast})^{-\frac{1}{2}}A(x-\mu^{k-1})]_{j}|\\
&+\max_{j\in S^{c}}|[(A^{\ast}(AA^{\ast})^{-1}e]_{j}|\\
&\leq\max_{j\in S^{c}}|\langle((AA^{\ast})^{-\frac{1}{2}}A)_{j},(AA^{\ast})^{-\frac{1}{2}}A(x-\mu^{k-1})\rangle|+\sqrt{1+\theta}\|e\|_{2}.
\end{aligned}
\end{equation}
Combining (\ref{equa1}) and (\ref{equa2}) gives
\begin{equation}\label{equa9}
\begin{aligned}
&\overline{\mathbb{P}}(\phi_{k}\geq\psi_{k})\\
&\leq\overline{\mathbb{P}}[\max_{j\in S^{c}}|\langle((AA^{\ast})^{-\frac{1}{2}}A)_{j},(AA^{\ast})^{-\frac{1}{2}}A(x-\mu^{k-1})\rangle|\\
&\geq\frac{1}{\sqrt{s-k+1}}(\|x_{T_{s-k+1}}\|_{2}-\gamma\|x-\mu^{k-1}\|_{2})\\
&-(\sqrt{1+\theta}/\sqrt{s-k+1}+\sqrt{1+\theta})\|e\|_{2}].
\end{aligned}
\end{equation}
Since $\|x_{T_{s-k+1}}\|\geq\sqrt{s-k+1}\widehat{x}_{s}$, $\|e\|_{2}\leq\epsilon\widehat{x}_{s}$,
\begin{equation*}
\begin{aligned}
&\frac{1}{\sqrt{s-k+1}}\|x_{T_{s-k+1}}\|_{2}-(\sqrt{1+\theta}/\sqrt{s-k+1}+\sqrt{1+\theta})\|e\|_{2}\\
&\geq\widehat{x}_{s}-\epsilon(\sqrt{1+\theta}/\sqrt{s-k+1}+\sqrt{1+\theta})\widehat{x}_{s}.\\
\end{aligned}
\end{equation*}
By {\em Lemma \ref{lemm6}} and since $\widehat{x}_{1}\leq\sigma\widehat{x}_{s}$, we have
\begin{equation*}
\begin{aligned}
\|x-\mu^{k-1}\|_{2}&\leq\frac{\|x_{T_{k-1}^{c}}\|_{2}}{\sqrt{1-\delta^{2}}}+\frac{\sqrt{1+\delta}\|e\|_{2}}{1-\delta}\\
&\leq\frac{\sqrt{s-k+1}}{\sqrt{1-\delta^{2}}}\widehat{x}_{1}+\frac{\sqrt{1+\delta}\|e\|_{2}}{1-\delta}\\
&\leq\frac{\sqrt{s-k+1}}{\sqrt{1-\delta^{2}}}\sigma\widehat{x}_{s}+\frac{\sqrt{1+\delta}}{1-\delta}\epsilon\widehat{x}_{s}.\\
\end{aligned}
\end{equation*}
With $\epsilon$, $\theta$ and $\gamma$ small enough, it can be derived
\begin{equation*}
\begin{aligned}
&1-\epsilon(\sqrt{1+\theta}/\sqrt{s-k+1}+\sqrt{1+\theta})\\
&\geq2\gamma(\frac{\sigma}{1-\delta^{2}}+\frac{\epsilon}{\sqrt{s-k+1}\sqrt{1-\delta}}).\\
\end{aligned}
\end{equation*}
Furthermore,
\begin{equation*}
\begin{aligned}
&\frac{1}{\sqrt{s-k+1}}\|x_{T_{s-k+1}}\|_{2}-(\sqrt{1+\theta}/\sqrt{s-k+1}+\sqrt{1+\theta})\|e\|_{2}\\
&\geq2\frac{\gamma}{\sqrt{s-k+1}}\|x-\mu_{k-1}\|_{2}.\\
\end{aligned}
\end{equation*}
Together with (\ref{equa9}),
\begin{equation*}
\begin{aligned}
&\overline{\mathbb{P}}(\phi_{k}\geq\psi_{k})\\
&\leq\overline{\mathbb{P}}[\max_{j\in S^{c}}|\langle((AA^{\ast})^{-\frac{1}{2}}A)_{j},(AA^{\ast})^{-\frac{1}{2}}A(x-\mu^{k-1})\rangle|\\
&\geq\frac{\gamma}{\sqrt{s-k+1}}\|x-\mu_{k-1}\|_{2}]\\
&\leq\overline{\mathbb{P}}[\max_{j\in S^{c}}|\langle((AA^{\ast})^{-\frac{1}{2}}A)_{j},(AA^{\ast})^{-\frac{1}{2}}A(x-\mu^{k-1})\rangle|\\
&\geq\frac{\gamma}{\sqrt{s}}\|x-\mu_{k-1}\|_{2}]\\
\end{aligned}
\end{equation*}
Exploiting the fact $\|(AA^{\ast})^{-\frac{1}{2}}A(x-\mu^{k-1})\|_{2}\leq\sqrt{1+\gamma}\|(x-\mu^{k-1})\|_{2}$, we have
\begin{equation*}
\begin{aligned}
&\overline{\mathbb{P}}(\phi_{k}\geq\psi_{k})\\
&\leq\overline{\mathbb{P}}(\max_{j\in S^{c}}|\langle((AA^{\ast})^{-\frac{1}{2}}A)_{j},(AA^{\ast})^{-\frac{1}{2}}A(x-\mu^{k-1})\|_{2}\\
&\geq\frac{\gamma}{\sqrt{1+\gamma}\sqrt{s}}\|(AA^{\ast})^{-\frac{1}{2}}A(x-\mu^{k-1})\|_{2})\\
&\leq\sum_{j\in S^{c}}\mathbb{P}(|\langle((AA^{\ast})^{-\frac{1}{2}}A)_{j},(AA^{\ast})^{-\frac{1}{2}}A(x-\mu^{k-1})\|_{2}\\
&\geq\frac{\gamma}{\sqrt{1+\gamma}\sqrt{s}}\|(AA^{\ast})^{-\frac{1}{2}}A(x-\mu^{k-1})\|_{2})\\
&\leq4(N-s)\exp(-\frac{C_{4}'\gamma^{2}M}{(1+\gamma)s}).
\end{aligned}
\end{equation*}
The last step is due to {\em Remark \ref{remark}}. Consequently, the failure probability $P$ satisfies
\begin{equation*}
\begin{aligned}
P&\leq6(N-s)\exp(-\min\{C_{2},C_{2}',C_{2}''\}\min\{\delta^{2},\gamma^{2},\theta^{2}\}M))\\
&+4(N-s)\exp(-\frac{C_{4}'\gamma^{2}M}{(1+\gamma)s})\\
&\leq6(N-1)\exp(-\min\{C_{2},C_{2}',C_{2}''\}\min\{\delta^{2},\gamma^{2},\theta^{2}\}M))\\
&+N^{2}\exp(-\frac{C_{4}'\gamma^{2}M}{(1+\gamma)s}).\\
\end{aligned}
\end{equation*}
With failure probability at most $P$, $\psi_{s}>\phi_{s}$ hold. As a result, $T_{s}=S$ and $\mu^{s}$ is a feasible solution. By {\em Lemma \ref{Pro1}}, the inequality
\begin{equation*}
\begin{aligned}
\|x-\mu^{s}\|_{2}&\leq\frac{1}{1-\delta}\|A(x-\mu^{s})\|_{2}\\
&\leq\frac{1}{\sqrt{1-\delta}}(\|y-Ax\|_{2}+\|y-A\mu^{s}\|_{2})\\
&\leq\frac{1}{\sqrt{1-\delta}}\|e\|_{2}
\end{aligned}
\end{equation*}
fails with probability at most $2\exp(-C_{2}\delta^{2}M)$. The final failure probability can be bounded by
\begin{equation*}
\begin{aligned}
&6N\exp(-\min\{C_{2},C_{2}',C_{2}''\}\min\{\delta^{2},\gamma^{2},\theta^{2}\}M))\\
&+N^{2}\exp(-\frac{C_{4}'\gamma^{2}M}{(1+\gamma)s})\\
&\leq 6N^{2}\exp(-\frac{C'M}{s})\\
&\leq6N^{-\alpha}.\\
\end{aligned}
\end{equation*}
The last step holds when $M\geq\frac{2+\alpha}{C'}s\ln(N)$.

 In the previous proof, the RIPs of $A$, $(AA^{\ast})^{-\frac{1}{2}}A$ and $(AA^{\ast})^{-1}A$ are considered individually. It is possible to exploit the relation of these properties. The true probabilities are larger than the ones in {\em Proposition \ref{theorem1}}, though it is not a focus of this article.
\section{Conclusions}
The NST+HT+FB and AdptNST+HT+FB algorithms are designed to find sparse solutions of under-determined linear systems. The convergence result of NST+HT+FB and numerical experiments about the effectiveness and the speed of NST+HT+FB have been presented in \cite{R25}. In this paper, the theoretical analysis of convergence results for both NST+HT+FB and AdptNST+HT+FB has been further elaborated. Our analysis improves the RIP and P-RIP condition of NST+HT+FB from $\delta_{2s}+\sqrt{2}\gamma_{3s}<1$ to $\delta_{2s}^{2}+2\gamma_{3s}^{2}<1$ and demonstrates that AdptNST+HT+FB converges in finitely many steps. The number of iterations for recovering an $s$-sparse signal of the two algorithms are also derived. In addition, we show that the number of iterations can be significantly lowered by exploiting the structure of the specific sparse signal or the random matrix.

\ack This work was partially supported by the NSF of USA (DMS-1313490, DMS-1615288), the China Scholarship Council, the National Natural Science Foundation of China (Grant Nos.61379014) and the Natural Science Foundation of Tianjin ( No.16JCYBJC15900).


\begin{thebibliography}{9}

\bibitem{R1} Li S, Qian T. On sparse representation of analytic signal in Hardy space[J]. Mathematical Methods in the Applied Sciences, 2013, 36(17): 2297-2310.
\bibitem{R2}    Haupt J, Bajwa W U, Rabbat M, et al. Compressed sensing for networked data[J]. IEEE Signal Processing Magazine, 2008, 25(2): 92-101.
\bibitem{R3}    Parvaresh F, Vikalo H, Misra S, et al. Recovering sparse signals using sparse measurement matrices in compressed DNA microarrays[J]. IEEE Journal of Selected Topics in Signal Processing, 2008, 2(3): 275-285.
\bibitem{R4}   Gedalyahu K, Eldar Y C. Time-delay estimation from low-rate samples: A union of subspaces approach[J]. IEEE Transactions on Signal Processing, 2010, 58(6): 3017-3031.
\bibitem{R5}   Natarajan B K. Sparse approximate solutions to linear systems[J]. SIAM journal on computing, 1995, 24(2): 227-234.
\bibitem{R6}   Chen S S, Donoho D L, Saunders M A. Atomic decomposition by basis pursuit[J]. SIAM review, 2001, 43(1): 129-159.
\bibitem{R7}   Tropp J A. Just relax: Convex programming methods for identifying sparse signals in noise[J]. IEEE transactions on information theory, 2006, 52(3): 1030-1051.
\bibitem{R8}    Mi W, Qian T, Li S. Basis pursuit for frequency domain identification[J]. Mathematical Methods in the Applied Sciences, 2016, 39(3): 498-507.
\bibitem{R9}   Chartrand R, Staneva V. Restricted isometry properties and nonconvex compressive sensing[J]. Inverse Problems, 2008, 24(3): 035020.
\bibitem{R10}   Foucart S, Lai M J. Sparsest solutions of underdetermined linear systems via $\ell_{q}$-minimization for $0<q\leq1$[J]. Applied and Computational Harmonic Analysis, 2009, 26(3): 395-407.
\bibitem{R11}  Pant J K, Lu W S, Antoniou A. New Improved Algorithms for Compressive Sensing Based on $\ell_{p}$ Norm[J]. IEEE Transactions on Circuits and Systems II: Express Briefs, 2014, 61(3): 198-202.
 \bibitem{R12}   Ambat S K, Hari K V S. An iterative framework for sparse signal reconstruction algorithms[J]. Signal Processing, 2015, 108: 351-364.

 \bibitem{R13}   Tibshirani R. Regression shrinkage and selection via the lasso[J]. Journal of the Royal Statistical Society. Series B (Methodological), 1996: 267-288.
 \bibitem{R14}   Gorodnitsky I F, Rao B D. Sparse signal reconstruction from limited data using FOCUSS: A re-weighted minimum norm algorithm[J]. IEEE Transactions on signal processing, 1997, 45(3): 600-616.
 \bibitem{R15}   Yin W, Osher S, Goldfarb D, et al. Bregman iterative algorithms for $\ell_{1}$-minimization with applications to compressed sensing[J]. SIAM Journal on Imaging sciences, 2008, 1(1): 143-168.
\bibitem{R16}    Mallat S G, Zhang Z. Matching pursuits with time-frequency dictionaries[J]. IEEE Transactions on signal processing, 1993, 41(12): 3397-3415.
 \bibitem{R17}   Pati Y C, Rezaiifar R, Krishnaprasad P S. Orthogonal matching pursuit: Recursive function approximation with applications to wavelet decomposition[C]Signals, Systems and Computers, 1993. 1993 Conference Record of The Twenty-Seventh Asilomar Conference on. IEEE, 1993: 40-44.
 \bibitem{R18}   Donoho D L, Tsaig Y, Drori I, et al. Sparse solution of underdetermined systems of linear equations by stagewise orthogonal matching pursuit[J]. IEEE Transactions on Information Theory, 2012, 58(2): 1094-1121.
 \bibitem{R19}    Tropp J, Needell D, Vershynin R. Iterative signal recovery from incomplete and inaccurate measurements[C]Proc. Information Theory and Applications Workshop. 2008.
 \bibitem{R20}      Dai W, Milenkovic O. Subspace pursuit for compressive sensing signal reconstruction[J]. IEEE Transactions on Information Theory, 2009, 55(5): 2230-2249.
 \bibitem{R21}   Blumensath T, Davies M E. Iterative hard thresholding for compressed sensing[J]. Applied and computational harmonic analysis, 2009, 27(3): 265-274.
 \bibitem{R22}    Daubechies I, Defrise M, De Mol C. An iterative thresholding algorithm for linear inverse problems with a sparsity constraint[J]. Communications on pure and applied mathematics, 2004, 57(11): 1413-1457.
 \bibitem{R23}     Foucart S. Hard thresholding pursuit: an algorithm for compressive sensing[J]. SIAM Journal on Numerical Analysis, 2011, 49(6): 2543-2563.
 \bibitem{R24}   Bouchot J L, Foucart S, Hitczenko P. Hard thresholding pursuit algorithms: number of iterations[J]. Applied and Computational Harmonic Analysis, 2016, 41(2): 412-435.
 \bibitem{R25}   Li S, Liu Y, Mi T. Fast thresholding algorithms with feedbacks for sparse signal recovery[J]. Applied and Computational Harmonic Analysis, 2014, 37(1): 69-88.
 \bibitem{R26}  Candes E J, Tao T. Decoding by linear programming[J]. IEEE transactions on information theory, 2005, 51(12): 4203-4215.
 \bibitem{R27}  Lin J, Li S. Nonuniform support recovery from noisy random measurements by orthogonal matching pursuit[J]. Journal of Approximation Theory, 2013, 165(1): 20-40.
 \bibitem{R28}    Tropp J A, Gilbert A C. Signal recovery from random measurements via orthogonal matching pursuit[J]. IEEE Transactions on information theory, 2007, 53(12): 4655-4666.
 \bibitem{R29}     Foucart S, Rauhut H. A mathematical introduction to compressive sensing[M]. Basel: Birkhauser, 2013.

\end{thebibliography}
\end{document}